\documentclass[MSNbibl,nameyear,citesort,seceqn,dvips]{arxbj}
\usepackage{graphicx}
\aid{0}
\volume{19}
\issue{4}
\pubyear{2013}
\firstpage{1176}
\lastpage{1211}
\doi{10.3150/13-BEJSP16} 

\makeatletter
\define@key{bid}{pmid}{}
\define@key{bid}{pmcid}{}

\newcommand{\indep}{\perp\hspace*{-6.2pt}\perp}

\newcommand{\dd}{\mathrm{d}}
\newcommand{\Expo}{\operatorname{Expo}}
\newcommand{\Poisson}{\operatorname{Poisson}}
\newcommand{\Var}{\operatorname{Var}}
\newcommand{\sign}{\operatorname{sign}}

\newcommand{\given}{ | }

\newtheorem{theorem}{Theorem}
\newproclaim{definition}{Definition}

\newcommand{\param}{\theta}
\newcommand{\pobs}{p_Y}
\newcommand{\psci}{p_X}
\newcommand{\pproc}{p_T}
\newcommand{\pwork}{\tilde{p}_X}
\newcommand{\pywork}{\tilde{p}_Y}
\newcommand{\obsparam}{\xi}
\newcommand{\workparam}{\eta}

\newcommand{\est}{\hat{\param}} 

\newcommand{\proc}{\mathcal{P}}
\makeatother

\begin{document}
\begin{frontmatter}

\title{The potential and perils of preprocessing: Building new foundations}

\runtitle{Potential and perils of preprocessing}

\begin{aug}
\author[a]{\fnms{Alexander W.} \snm{Blocker}\thanksref{e1}\ead[label=e1,mark]{ablocker@post.harvard.edu}}
\and
\author[a]{\fnms{Xiao-Li} \snm{Meng}\corref{}\thanksref{e2}\ead[label=e2,mark]{meng@stat.harvard.edu}}
\address[a]{Department of Statistics, Harvard University, Cambridge,
MA 02138,
USA.\\
\printead{e1,e2}}
\runauthor{A.W. Blocker and X.-L. Meng} 
\end{aug}


%
\begin{abstract}
Preprocessing forms an oft-neglected foundation for a wide range of
statistical and scientific analyses.
However, it is rife with subtleties and pitfalls.
Decisions made in preprocessing constrain all later analyses and are
typically irreversible.
Hence, data analysis becomes a collaborative endeavor by all parties
involved in data collection, preprocessing and curation, and downstream
inference.
Even if each party has done its best given the information and
resources available to them, the final result may still fall short of
the best possible in the traditional single-phase inference framework.
This is particularly relevant as we enter the era of ``big data''.
The technologies driving this data explosion are subject to complex new
forms of measurement error.
Simultaneously, we are accumulating increasingly massive databases of
scientific analyses.
As a result, preprocessing has become more vital (and potentially more
dangerous) than ever before.

We propose a theoretical framework for the analysis of preprocessing
under the banner of multiphase inference.
We provide some initial theoretical foundations for this area,
including distributed preprocessing, building upon previous work in
multiple imputation.
We motivate this foundation with two problems from biology and
astrophysics, illustrating multiphase pitfalls and potential solutions.
These examples also emphasize the motivations behind multiphase
analyses---both practical and theoretical.
We demonstrate that multiphase inferences can, in some cases, even
surpass standard single-phase estimators in efficiency and robustness.
Our work suggests several rich paths for further research into the
statistical principles underlying preprocessing.
%
%
To tackle our increasingly complex and massive data, we must ensure
that our inferences are built upon solid inputs and sound principles.
Principled investigation of preprocessing is thus a vital direction for
statistical research.
\end{abstract}

%
\begin{keyword}
\kwd{data compression}
\kwd{data repositories}
\kwd{measurement error}
\kwd{multiphase inference}
\kwd{multiple imputation}
\kwd{statistical principles}
\end{keyword}

\end{frontmatter}

\section{What is multiphase inference?}

\subsection{Defining multiphase problems}
Preprocessing and the analysis of preprocessed data are ubiquitous
components of statistical inference, but their treatment has often been informal.
We aim to develop a theory that provides a set of formal statistical
principles for such problems under the banner of multiphase inference.
The term ``multiphase'' refers to settings in which inferences are
obtained through the application of multiple procedures in sequence,
with each procedure taking the output of the previous phase as its input.
This encompasses settings such as multiple imputation (MI, \citet{Rubin1987}) and extends to other situations.
In a multiphase setting, information can be passed between phases in an
arbitrary form; it need not consist of (independent) draws from a
posterior predictive distribution, as is typical with multiple imputation.
Moreover, the analysis procedure for subsequent phases is not
constrained to a particular recipe, such as Rubin's MI combining rules
(\citeyear{Rubin1987}).

The practice of multiphase inference is currently widespread in applied
statistics.
It is widely used as an analysis technique within many
publications---any paper that uses a ``pipeline'' to obtain its final
inputs or clusters estimates from a previous analysis provides an example.
Furthermore, projects in astronomy, biology, ecology, and social
sciences (to name a small sampling) increasingly focus on building
databases for future analyses as a primary objective.
These projects must decide what levels of preprocessing to apply to
their data and what additional information to provide to their users.
Providing all of the original data clearly allows the most flexibility
in subsequent analyses.
In practice, the journey from raw data to a complete model is typically
too intricate and problematic for the majority of users, who instead
choose to use preprocessed output.

Unfortunately, decisions made at this stage can be quite treacherous.
Preprocessing is typically irreversible, necessitating assumptions
about both the observation mechanisms and future analyses. These
assumptions constrain all subsequent analyses.
Consequently, improper processing can cause a disproportionate amount
of damage to a whole body of statistical results.
However, preprocessing can be a powerful tool.
It alleviates complexity for downstream researchers, allowing them to
deal with smaller inputs and (hopefully) less intricate models.
This can provide large mental and computational savings.

Two examples of such trade-offs come from NASA and high-throughput biology.
When NASA satellites collect readings, the raw data are usually massive.
These raw data are referred to as the ``Level 0'' data (\citet{ChandraSDP}).
The Level 0 data are rarely used directly for scientific analyses.
Instead, they are processed to Levels 1, 2, and 3, each of which
involves a greater degree of reduction and adjustment.
Level 2 is typically the point at which the processing becomes irreversible.
\citet{Braverman2012} provide an excellent illustration of this process
for the Atmospheric Infrared Sounder (AIRS) experiment.
This processing can be quite controversial within the astronomical community.
Several upcoming projects, such as the Advanced Technology Solar
Telescope (ATST) will not be able to retain the Level 0 or Level 1 data
(\citet{Davey2012}).
This inability to obtain raw data and increased dependence on
preprocessing has transformed low-level technical issues of calibration
and reduction into a pressing concern.

High-throughput biology faces similar challenges.
Whereas reproducibility is much needed (e.g., \citet{Ioannidis2011}),
sharing raw datasets is difficult because of their sizes.
The situation within each analysis is similar.
Confronted with an overwhelming onslaught of raw data, extensive
preprocessing has become crucial and ubiquitous.
Complex models for genomic, proteomic, and transcriptomic data are
usually built upon these heavily-processed inputs.
This has made the intricate details of observation models and the
corresponding preprocessing steps the groundwork for entire fields.

To many statisticians, this setting presents something of a conundrum.
After all, the ideal inference and prediction will generally use a
complete correctly-specified model encompassing the underlying process
of interest and all observation processes.
Then, why are we interested in multiphase?
We focus on settings where there is a natural separation of knowledge
between analysts, which translates into a separation of effort.
The first analyst(s) involved in preprocessing often have better
knowledge of the observation model than those performing subsequent analyses.
For example, the first analyst may have detailed knowledge of the
structure of experimental errors, the equipment used, or the
particulars of various protocols.
This knowledge may not be easy to encapsulate for later analysts---the
relevant information may be too large or complex, or the methods
required to exploit this information in subsequent analyses may be
prohibitively intricate.
%
%
Hence, the practical objective in such settings is to enable the best
possible inference given the constraints imposed and provide an account
of the trade-offs and dangers involved.
To borrow the phrasing of \citet{Meng2003} and \citet{Rubin1996}, we
aim for achievable practical efficiency rather than theoretical
efficiency that is practically unattainable.

Multiphase inference currently represents a serious gap between
statistical theory and practice.
We typically delineate between the informal work of preprocessing and
feature engineering and formal, theoretically-motivated work of
estimation, testing, and so forth.
However, the former fundamentally constrains what the latter can accomplish.
As a result, we believe that it represents a great challenge and
opportunity to build new statistical foundations to inform statistical practice.

\subsection{Practical motivations}
\label{sec:examples}

We present two examples that show both the impetus for and perils of
undertaking multiphase analyses in place of inference with a complete,
joint model. The first concerns microarrays, which allow the analysis
of thousands of genes in parallel.
We focus on expression microarrays, which measure the level of gene
expression in populations of cells based upon the concentration of RNA
from different genes.
These are typically used to study changes in gene expression between
different experimental conditions.

In such studies, the estimand of interest is typically the log-fold
change in gene expression between conditions.
However, the raw data consist only of intensity measurements for each
probe on the array, which are grouped by gene along with some form of controls.
These intensities are subject to several forms of observation noise,
including additive background variation and additional forms of
interprobe and interchip variation (typically modeled as multiplicative noise).
To deal with these forms of observation noise, a wide range of
background correction and normalization strategies have been developed
(for a sampling,
see \citet{Tusher2001},
\citet{Quackenbush2002},
\citet{Affymetrix2002},
\citet{Irizarry2003},
\citet{McGee2006},
\citet{Ritchie2007},
\citet{Xie2009}).
Later analyses then focus on the scientific question of interest
without, for the most part, addressing the underlying details of the
observation mechanisms.

Background correction is a particularly crucial step in this process,
as it is typically the point at which the analysis moves from the
original intensity scale to the log-transformed scale.
As a result, it can have a large effect on subsequent inferences about
log-fold changes, especially for genes with low expression levels in
one condition (\citet{Smyth2005},
\citet{Irizarry2006}).
One common method (MAS5), provided by one microarray manufacturer, uses
a combination of background subtraction and truncation at a fixed lower
threshold for this task (\citet{Affymetrix2002}). Other more
sophisticated techniques use explicit probability models for this
de-convolution.
A model with normally-distributed background variation and
exponentially distributed expression levels has proven to be the most
popular in this field (\citet{McGee2006},
\citet{Xie2009}).

Unfortunately, even the most sophisticated available techniques pass
only point estimates onto downstream analyses.
This necessitates ad-hoc screening and corrections in subsequent
analyses, especially when searching for significant changes in
expression (e.g., \citet{Tusher2001}).
Retaining more information from the preprocessing phases of these
analyses would allow for better, simpler inference techniques with
greater power and fewer hacks.
The motivation behind the current approach is quite understandable:
scientific investigators want to focus on their processes of interest
without becoming entangled in the low-level details of observation mechanisms.
Nevertheless, this separation can clearly compromise the validity of
their results.

The role of preprocessing in microarray studies extends well beyond
background correction.
Normalization of expression levels across arrays, screening for data
corruption, and other transformations preceding formal analysis are standard.
Each technique can dramatically affect downstream analyses.
For instance, quantile normalization equates quantiles of expression
distributions between arrays, removing a considerable amount of information.
This mutes systematic errors (\citet{Bolstad2003}), but it can seriously
compromise analyses in certain contexts (e.g., miRNA studies).

Another example of multiphase inference can be found in the estimation
of correlations based upon indirect measurements.
This appears in many fields, but astrophysics provides one recent and
striking case.
The relationships between the dust's density, spectral properties, and
temperature are of interest in studies of star-forming dust clouds.
These characteristics shed light on the mechanisms underlying star
formation and other astronomical processes.
Several studies (e.g., \citet{Dupac2003},
\citet{Desert2008},
\citet{Anderson2010},
\citet{Paradis2010}) have
investigated these relationships, finding negative correlations between
the dust's temperature and spectral index.
This finding is counter to previous astrophysical theory, but it has
generated many alternative explanations.

Such investigations may, however, be chasing a phantasm.
These correlations have been estimated by simply correlating point
estimates of the relevant quantities (temperature $T$ and spectral
index $\beta$) based on a single set of underlying observations.
As a result, they may conflate properties of this estimation procedure
with the underlying physical mechanisms of interest.
This has been noted in the field by \citet{Shetty2009}, but the
scientific debate on this topic continues.
\citet{Kelly2012} provide a particularly strong argument, using a
cohesive hierarchical Bayesian approach, that improper multiphase
analyses have been a pervasive issue in this setting.
Improper preprocessing led to incorrect, negative estimates of the
correlation between temperature and spectral index, according to \citet
{Kelly2012}.
These incorrect estimates even appeared statistically significant with
narrow confidence intervals based on standard methods.
On a broader level, this case again demonstrates some of the dangers of
multiphase analyses when they are not carried out properly.
Those analyzing this data followed an intuitive strategy: estimate what
we want to work with ($T$ and $\beta$), then use it to estimate the
relationship of interest.
Unfortunately, such intuition is not a recipe for valid statistical
inference.

\subsection{Related work}

Multiphase inference has wide-ranging connections to both the
theoretical and applied literatures.
It is intimately related to previous work on multiple imputation and
missing data (\citeauthor{Rubin1976} (\citeyear{Rubin1976,Rubin1987,Rubin1996}),
\citet{Meng1994},
\citet{Meng2003},
\citet{Xie2012}).
In general, the problem of multiphase inference can be formulated as
one of missing data.
However, in the multiphase setting, missingness arises from the
preprocessing choices made, not a probabilistic response mechanism.
Thus, we can leverage the mathematical and computational methods of
this literature, but many of its conceptual tools need to be modified.
Multiple imputation addresses many of the same issues as multiphase
inference and is indeed a special case of the latter.
Concepts such as congeniality between imputation and analysis models
and self-efficiency (\citet{Meng1994}) have natural analogues and roles
to play in the analysis of multiphase inference problems.

Multiphase inference is also tightly connected to work on the
comparison of experiments and approximate sufficiency, going back to
\citeauthor{Blackwell51} (\citeyear{Blackwell51,Blackwell53}) and continuing through \citet
{LeCam1964} and \citet{Goel79}, among others.
This literature has addressed the relationship between decision
properties and the probabilistic structure of experiments, the
relationship between different notions of statistical information, and
notions of approximate sufficiency---all of these are quite relevant
for the study of multiphase inference.
We view the multiphase setting as an extension of this work to address
a broader range of real-world problems, as we will discuss in
Section~\ref{sec:riskmonotone}.

The literature on Bayesian combinations of experts also informs our
thinking on multiphase procedures.
\citet{Kadane1993} provides an excellent review of the field, while
\citet{Lindley1979} provides the core formalisms of interest for the
multiphase setting.
Overall, this literature has focused on obtaining coherent (or
otherwise favorable) decision rules when combining information from
multiple Bayesian agents, in the form of multiple posterior
distributions. We view this as a best-case scenario, focusing our
theoretical development towards the mechanics of passing information
between phases.
We also focus on the sequential nature of multiphase settings and the
challenges this brings for both preprocessors and downstream analysts,
in contrast to the more ``parallel'' or simultaneous focus of the
literature mentioned above.

There are also fascinating links between multiphase inference and the
signal processing literature.
There has been extensive research on the design of quantizers and other
compression systems; see for example \citet{Gray1998}.
Such work is often focused on practical questions, but it has also
yielded some remarkable theory.
In particular, the work of \citet{Nguyen2009} on the relationship
between surrogate loss functions in quantizer design and
$f$-divergences suggests possible ways to develop and analyze a wide
class of multiphase procedures, as we shall discuss in Section~\ref{sec:future}.

\section{Multiphase logic and concepts for preprocessing}
\label{sec:concepts}

\subsection{A model for two phases}
\label{sec:model}

To formalize the notion of multiphase inference, we begin with a formal
model for two-phase settings.
The first phase consists of the data generation, collection, and
preprocessing, while the second phase consists of inference using the
output from the first phase.
We will call the first-phase agent the ``preprocessor'' and the
second-phase agent the ``downstream analyst''.
The preprocessor observes the raw data $Y$.
This is a noisy realization of $X$, variables of interest that are not
directly obtainable from a given experiment, e.g., gene expression from
sequencing data, or stellar intensity from telescopic observations.

We assume that the joint density of $X$ and $Y$ with respect to product
measure $\mu_X \times\mu_Y$ can be factored as
%
%
\begin{equation}
\label{e:model} p_{Y,X}(Y,X \given\param, \obsparam) = \pobs(Y \given X,
\param, \obsparam) \cdot\psci(X \given\param, \obsparam) = \pobs(Y \given X,
\obsparam) \cdot\psci(X \given\param)  .
\end{equation}
Here, $\psci$ encapsulates the underlying process of interest and
$\pobs
$ encapsulates the observation process.
We assume that $\param$ is of fixed dimension in all asymptotic settings.
In practice, the preprocessor should be able to postulate a reasonable
``observation model'' $\pobs(Y \given X, \obsparam)$, but will not
always know the true ``scientific model'' $\psci(X \given\param)$.
This is analogous to the MI setting, where the imputer does not know
the form of the final analysis.

\begin{figure}[b]

\includegraphics{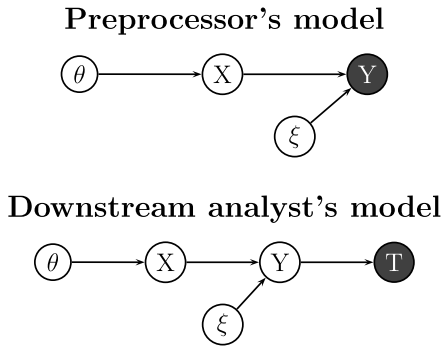}

\caption{Graphical diagram of our generic two-phase setting. The
preprocessor observes $Y$ from the original data
generating process and outputs $T$, with $X$ as missing data. The
downstream analyst observes the preprocessor's output $T$ and has
both $X$ and $Y$ missing.}\label{fig:models}
\end{figure}

Using this model, the preprocessor provides the downstream analyst with
some output $T = T(Y, U)$, where $U$ is a (possibly stochastic)
additional input.
When $T(Y, U)$ is stochastic (e.g., an MCMC output), the conditional
distribution $\pproc(T|Y)$ is its theoretical description instead of
its functional form.
However, for simplicity, we will present our results when $T$ is a
deterministic function of $Y$ only, but many results generalize easily.
Given such $T$, downstream analysts can carry out their inference procedures.
Figure~\ref{fig:models} depicts our general model setup.

This model incorporates several restrictions.
First, it is Markovian with respect to $Y$, $X$, and $\param$; $Y$ is
conditionally independent of $\param$ given $X$ (and $\obsparam$).
Second, the parameters governing the observation process ($\obsparam$)
and those governing the scientific process ($\param$) are distinct.
In Bayesian settings, we further assume that $\obsparam$ and $\param$
are independent \textit{a priori}.
The parameters $\obsparam$ are nuisance from the perspective of all
involved; the downstream analyst wants to draw inferences about $X$ and
$\param$, and the preprocessor wants to pass forward information that
will be useful for said inferences.
If downstream inferences are Bayesian with respect to $\obsparam$, then
$\pobs(Y|X) = \int\pobs(Y | X, \obsparam) \pi_\obsparam(\obsparam
)\,\dd
\mu_{\obsparam}(\obsparam)$ (which holds under~(\ref{e:model})) is
sufficient for all inference under the given model and prior.
Hence, this conditional density is frequently of interest in our
theoretical development, as is the corresponding marginalized model
$p_{X,Y}(Y, X | \param) = \int\pobs(Y | X, \obsparam) \psci(X |
\param
) \pi_\obsparam(\obsparam)\,\dd\mu_{\obsparam}(\obsparam)$.
We will compare results obtained with a fixed prior to those obtained
in a more general setting to better understand the effects of nuisance
parameters in multiphase inference.

These restrictions are somewhat similar to those underlying Rubin's
(\citeyear{Rubin1976}) definition of ``missing at random''; however, we
do not have missing data mechanism (MDM) in this setting \emph{per se}.
The distinction between missing and observed data ($X$ and $Y$) is
fixed by the structure of our model.
In place of MDM, we have two imposed patterns of missingness: one for
the data-generating process, and one for the inference process.
The first is $\pobs(Y \given X, \obsparam)$, which creates a noisy
version of the desired scientific variables.
Here, $X$ can be considered the missing data and $Y$ the observed.
For the inference process, the downstream analyst observes $T$ in place
of $Y$ but desires inference for $\param$ based upon $\psci(X |
\param)$.
Hence, $Y$ and $X$ are both missing for the downstream analyst.
Neither pattern is entirely intrinsic to the problem---both are fixed
by choice.
The selection of scientific variables $X$ for a given marginal
likelihood $\pobs(Y | \param, \obsparam) = \int\pobs(Y | X,
\obsparam)
\psci(X | \param) \,\dd\mu_{X}(X)$ is a modeling decision.
The selection of preprocessing $T(Y)$ is a design decision.
This contrasts with the typical missing data setting, where MDM is
forced upon the analyst by nature.
With multiphase problems, we seek to design and evaluate engineered missingness.
Thus the investigation of multiphase inference requires tools and ideas
from design, inference, and computation in addition to the established
theory of missing data.

\subsection{Defining multiphase procedures}
\label{sec:def}

With this model in place, we turn to formally defining multiphase procedures.
This is more subtle than it initially appears.
In the MI setting, we focus on complete-data procedures for the
downstream analyst's estimation and do not restrict the dependence
structure between missing data and observations.
In contrast, we restrict the dependence structure as in (\ref
{e:model}), but place far fewer constraints on the analysts' procedures.
Here, we focus our definitions and discussion on the two-phase case of
a single preprocessor and downstream analyst.
This provides the formal structure to describe the interface between
any two phases in a chain of multiphase analyses.

In our multiphase setting, downstream analysts need not have any
complete-data procedure in the sense of one for inferring $\param$ from
$X$ and $Y$; indeed, they need not formally have one based only upon
$X$ for inferring $\param$.
We require only that they have a set of procedures for their desired
inference using the quantities provided from earlier phases as inputs
($T$), not necessarily using direct observations of $X$ or $Y$.
Such situations are common in practice, as methods are often built
around properties of preprocessed data such as smoothness or sparsity
that need not hold for the actual values of $X$.

For the preprocessor, the input is $Y$ and the output is $T$.
Here $T$ could consist of a vector of means with corresponding
standard errors, or, for discrete $Y$, $T$ could consist of carefully
selected cross-tabulations.
In general, $T$ clearly needs to be related to $X$ to capture
inferential information, but its actual form is influenced by practical
constraints (e.g., aggregation to lower than desired resolutions due to
data storage capacity).

For the downstream analyst, the input is $T$ and the output is an
inference for $\theta$.
This analyst can obviously adapt.
For example, suppose $\param= E(X_i)$ for each entry $i$ of $X$.
If the preprocessor provides $T_0 = \hat{X}$, the analyst may simply
use an unweighted mean to estimate $\param$.
If the preprocessor instead gives the analyst $T_1 = (\hat{X}, S)$,
where $S$ contains standard errors, the latter could instead use a
weighted mean to estimate $\param$.
This adaptation extends to an arbitrary number of possible inputs
$T_k$, each of which corresponds to a set of constraints facing the
preprocessor.

To formalize this notion of adaptation, we first define an index set
$C$ with one entry for each such set of constraints.
This maps between forms of input provided by the preprocessor and
estimators selected by the downstream analyst.
In this way, $C$ captures the downstream analyst's knowledge of
previous processing and the underlying probability model.
Thus, this index set plays an central role in the definition of
multiphase inference problems, far beyond that of a mere mathematical
formality; it regulates the amount of mutual knowledge shared between
the preprocessor and the downstream analyst.

Now, we turn to the estimators themselves.
We start with point estimation as a foundation for a broader class of problems.
Testing begins with estimating rejection regions, interval estimation
with estimating coverage, classification with estimating class
membership, and prediction with estimating future observations and,
frequently, intermediate parameters.
The framework we present therefore provides tools that can be adapted
for more than estimation theory.
We define multiphase estimation procedures as follows:
%
\begin{definition}
A \emph{multiphase estimation procedure} $\proc$ is a set of estimators
$\{\est_k(T_k) \dvt\break k \in C \}$ indexed by the set $C$, where $T_k$
corresponds to the output of the $k$th first-phase method; that is,
$\proc$ is a family of estimators with different inputs.
\end{definition}
When clear, we will drop the subscripts $k$ and index the estimators in
$\proc$ by their inputs.
This definition provides enough flexibility to capture many practical
issues with multiphase inference, and it can be iterated to define
procedures for analyses involving a longer sequence of preprocessors
and analysts.
It also encompasses the definition of a missing data procedure used by
\citet{Meng1994}.
Such procedures cannot, of course, be arbitrarily constructed if they
are to deliver results with general validity.
Hence, having defined these procedures, we will cull many of them from
consideration in Section~\ref{sec:riskmonotone}.

The obvious choice of our estimand, suggested by our notation thus far,
is the parameter for the scientific model, $\param$.
This is very amenable to mathematical analysis and relevant to many
investigations.
Hence, it forms the basis for our results in Section~\ref{sec:theory}.
However, for multiphase analyses, other classes of estimands may prove
more useful in practice.
In particular, functions of $X$, future scientific variables $X_{rep}$,
or future observations $Y_{rep}$ may be of interest.
Prediction of such quantities is a natural focus in the multiphase
setting because such statements are meaningful to both the preprocessor
and downstream analyst.
Such estimands naturally encompass a broad range of statistical
problems including prediction, classification, and clustering.
However, there is often a lack of mutual knowledge about $\psci(X
\given\param)$, so the preprocessor cannot expect to ``target''
estimation of $\param$ in general, as we shall discuss in Section~\ref{sec:remarks}.

\subsection{When is more better?}
\label{sec:riskmonotone}

It is not automatic for multiphase estimation procedures to produce
better results as the first phase provides more information.
To obtain a sensible context for theoretical development, we must
regulate the way that the downstream analyst adapts to different inputs.
For instance, they should obtain better results (in some sense) when
provided with higher-resolution information.
This carries over from the MI setting (\citet
{Meng1994}, \citet{Meng2003}, \citet{Meng2012},
\citet{Xie2012}), where notions such as
self-efficiency are useful for regulating the downstream analyst's procedures.
We define a similar property for multiphase estimation procedures, but
without restricting ourselves to the missing data setting.
Specifically, let $T_1 \preceq T_2$ indicate $T_1$ is a deterministic
function of $T_2$.
In practice, $T_1$ could be a subvector, aggregation, or other summary~of~$T_2$.

\begin{definition}[(Risk monotonicity)]
A multiphase estimation procedure $\proc$ is \emph{risk monotone} with
respect to a loss function $L$ if, for all pairs of outputs $T_1, T_2$,
$T_1 \preceq T_2$ implies $R(\est_2(T_2), L) \leq R(\est_1(T_1), L)$.
\end{definition}
An asymptotic analogue of risk monotonicity is defined as would be
expected, scaling the relevant risks at an appropriate rate to obtain
nontrivial limits.
This is a natural starting point for regulating multiphase estimation
procedures; stronger notions may be required for certain theoretical results.
Note that this definition does not require that ``higher-quality''
inputs necessarily lead to lower risk estimators.
Risk monotonicity requires only that estimators based upon a larger set
of inputs perform no worse than those with strictly less information
(in a deterministic sense).
However, risk monotonicity is actually quite tight in another sense.
It requires that additional information cannot be misused by the
downstream analyst, imposing a strong constraint on mutual knowledge.
For an example, consider the case of unweighted and weighted means.
To obtain better results when presented with standard errors, the
downstream analyst must know that they are being given (the correct)
standard errors and to weight by inverse variances.

This definition is related to the comparison of experiments, as
explored by \citeauthor{Blackwell51} (\citeyear{Blackwell51,Blackwell53}), but diverges on a
fundamental level.
Our ordering of experiments, based on deterministic functions, is more
stringent than that of \citet{Blackwell53}, but they are related.
Indeed, our $\preceq$ relation implies that of \citet{Blackwell53}.
In the latter work, an experiment $\alpha$ is defined as more
informative than experiment $\beta$, denoted $\alpha\supset\beta$, if
all losses attainable from $\beta$ are also attainable from $\alpha$.
This relation is also implied when $\alpha$ is sufficient for $\beta$.
Our stringency stems from our broader objectives in the multiphase setting.
From a decision-theoretic perspective, the partial ordering of
experiments investigated by Blackwell and others deal with which risks
are attainable given pairs of experiments, allowing for arbitrary
decision procedures.
In contrast, our criterion restricts procedures based on whether such
risks are actually attained, with respect to a particular loss function.
This is because, in the multiphase setting, it is not generally
realistic to expect downstream analysts to be capable of obtaining
optimal estimators for all forms of preprocessing.

The conceptually-simplest way to generate such a procedure is to begin
with a complete probability model for $\pobs(Y | \param)$.
Under traditional asymptotic regimes, all procedures consisting of
Bayes estimators based upon such a model will (with full knowledge of
the transformations involved in each $T_k$ and a fixed prior) be risk monotone.
The same is true asymptotically under the same regimes (for
squared-error loss) for procedures consisting of MLEs under a fixed model.
Under some other asymptotic regimes, however, these principles of
estimation do not guarantee risk-monotonicity; we explore this further
in Section~\ref{sec:missinfo}.
But such techniques are not the only way to generate risk monotone
procedures from probability models.
This is analogous to self-efficiency, which can be achieved by
procedures that are neither Bayesian nor MLE (\citet{Meng1994},
\citet{Xie2012}).

\begin{figure}[b]

\includegraphics{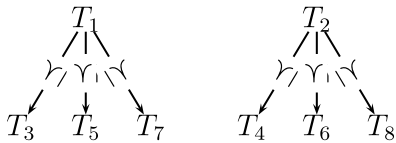}

\caption{Illustration of risk-monotone ``basis'' construction. In this
case, $T_1$ and $T_2$ form the basis set of
statistics. Each of these has three descendants ($T_3, T_5, T_7$ from
$T_1$ and $T_4, T_6, T_8$ from~$T_2$). These descendants are
deterministic functions of their parent, but they are not deterministic
functions of any other basis statistics. Given correctly-specified
models for $T_1$ and~$T_2$, a risk monotone procedure
can be constructed for all statistics ($T_1, \ldots, T_8$) shown here
as described in the text.}\label{fig:risk-monotone}
\end{figure}

A risk monotone procedure can be generated from any set of probability
models for distinct inputs that ``span'' the space of possible inputs.
Suppose that an analyst has a set of probability models, all correctly
specified, for $p_{T_b}(T_b | \param)$, where $b$ ranges over a subset
$B$ of the relevant index set $C$.
We also assume that this analyst has a prior distribution $\pi
_b(\param
)$ for each such basis models.
These priors need not agree between models; the analyst can build a
risk-monotone procedure from an inconsistent set of prior beliefs.
Suppose that the inputs $\{T_b : b \in B\}$ are not deterministic
functions of each other and all other inputs can be generated as
nontrivial deterministic transformations of one of these inputs.
Formally, we require $T_b \npreceq T_c$ for all distinct $b, c \in
B$ and, for each $k \in C$ there exists a unique $b \in B$ such that
$T_k \preceq T_b$ (each output is uniquely descended from a single $T_b$),
as illustrated in Figure~\ref{fig:risk-monotone}.
This set can form a basis, in a sense, for the given procedure.

Using the given probability models with a single loss function and set
of priors (potentially different for each model), the analyst can
derive a Bayes rule under each model.
For each $b \in B$, we require $\est(T_b)$ to be an appropriate Bayes
rule on said model.
As $T_k = g_k(T_b)$ for some function $g_k$, we then have the implied
$p_{T_k}(T_k | \param) = \int_{t : g_k(t) = T_k} p_{T_b}(t | \param)
\,\dd t$, yielding the Bayes rule for estimating $\param$ based on $T_k$,
which is no less risky than $\est(T_b)$.
The requirement that each output $T_k$ derives from a unique $T_b$
means that each basis component $T_b$ has a unique line of descendants.
Within each line, each descendant is comparable to only a single $T_b$
in the sense of deterministic dependence.
Between these lines, such comparisons are not possible.
This ensures the overall risk-monotonicity.

Biology provides an illustration of such bases.
A wide array of methodological approaches have been used to analyze
high-throughput gene expression data.
One approach, builds upon order and rank statistics (\citet
{Geman2004}, \citet{Geman2012}, \citet{Tan2005}).
Another common approach uses differences in gene expression between
conditions or experiments, often aggregating over pathways, replicates,
and so forth.
Each class of methods is based upon a different form of preprocessing:
ranks transformations for the former, normalization and aggregation for
the latter.
Taking procedures based on rank statistics and aggregate differences in
expression as a basis, we can consider constructing a risk-monotone
procedure as above.
Thus, the given formulation can bring together apparently disparate
methods as a first step in analyzing their multiphase properties.

Such constructions are, unfortunately, not sufficient to generate all
possible risk monotone procedures.
Obtaining more general conditions and constructions for risk monotone
procedures is a topic for further work.

\subsection{Revisiting our examples and probing our boundaries}
\label{sec:revisit}

By casting the examples in Section~\ref{sec:examples} into the formal
structure just established, we can clarify the practical role of each
mathematical component and see how to map theoretical results into
applied guidance.
We also provide an example that illustrates the boundaries of the
framework's utility, and another that demonstrates its formal limits.
These provide perspective on the trade-offs made in formalizing the
multiphase inference problem.

The case of microarray preprocessing presented previously fits quite
nicely into the model of Section~\ref{sec:model}.
There, $Y$ corresponds to the observed probe-level intensities,
$X$~corresponds to the true expression level for each gene under each
condition, and $\param$ corresponds to the parameters governing the
organism's patterns of gene expression.
In the microarray setting, $\pobs$ would characterize the relationship
between expression levels and observed intensities, governed by~$\obsparam$.
These nuisance parameters could include chip-level offsets, properties
of any additive background, and the magnitudes of other sources of variation.
The assumptions of a Markovian dependence structure and distinct
parameters for each part of the model appear quite reasonable in this
case, as (1) the observation $Y$ can only (physically) depend upon the
sample preparation, experimental protocol, and RNA concentrations in
the sample and (2) the distributions $\psci$ and $\pobs$ capture
physically distinct portions of the experiment.
Background correction, normalization, and the reduction of observations
to log-fold changes are common examples of preprocessing $T(Y)$. As
discussed previously, estimands based upon $X$ may be of greater
scientific interest than those based upon $\param$.
For instance, we may want to know whether gene expression changed
between two treatments in a particular experiment (a statement about~$X$)
than whether a parameter regulating the overall patterns of gene
expression takes on a particular value.

For the astrophysical example, the fit is similarly tidy.
The raw astronomical observations correspond to $Y$, the true
temperature, density, and spectral properties of each part of the dust
cloud become $X$, and the parameters governing the relationship between
these quantities (e.g., their correlation) form $\param$.
The $\pobs$ distribution governs the physical observation process,
controlled by~$\obsparam$.
This process typically includes the instruments' response to
astronomical signals, atmospheric distortions, and other earthbound phenomena.
As before, the conditional independence of $\param$ and $Y$ given $X$
and $\obsparam$ is sensible based upon the problem structure, as is the
separation of $\param$ and $\obsparam$. Here $X$ corresponds to signals
emitted billions or trillions of miles from Earth, whereas the
observation process occurs within ground- or space-based telescopes.
Hence, any non-Markovian effects are quite implausible.
Preprocessing $T(Y)$ corresponds to the (point) estimates of
temperature, density, and spectral properties from simple models of $Y$
given $X$ and~$\obsparam$.

The multiphase framework encompasses a broad range of settings, but it
does not shed additional light on all of them.
If $T$ is a many-to-one transformation of $Y$, then our framework
implies that the preprocessor and downstream analyst face structurally
different inference (and missing data) problems.
This is the essence of multiphase inference, in our view.
Settings where $\pobs(Y|X, \obsparam)$ is degenerate or $T$ is a
one-to-one function of $Y$ are boundary cases where our multiphase
interpretation and framework add little.

For a concrete example of these cases, consider a time-to-failure
experiment, with the times of failure $W_i \sim\mbox{i.i.d } \Expo
(\param)$, $i = 1, \ldots, n$.
Now, suppose that the experimenters actually ran the experiment in $m$
equally-sized batches.
They observe each batch only until its first failure; that is, they
observe and report $Y_b = \min\{W_i : i \mbox{ in batch } b \}$ for
each batch $b$.
Subsequent analysts have access only to $T = (Y_1, \ldots, Y_b)$.
This seems to be a case of preprocessing, but it actually resides at
the very edge of our framework.

We could take the complete observations to be $X$ and the batch minima
to be $Y$.
This would satisfy our Markov constraint, with a singular, and hence
deterministic, observation process $\pobs(Y | X)$ simply selecting a
particular order statistic within each batch.
However, $T(Y)$ is one-to-one; the preprocessor observes only the order
statistics, as does the downstream analyst.
There is no separation of inference between phases; the same quantities
are observed and missing to both the preprocessor and the downstream analyst.
Squeezing this case into the multiphase framework is technically valid
but unproductive.

The framework we present is not, however, completely generic.
Consider a chemical experiment involving a set of reactions.
The underlying parameters $\param$ describe the chemical properties
driving the reactions, $X$ are the actual states of the reaction, and
$Y$ are the (indirectly) measured outputs of the reactions.
The measurement process for these experiments, as described by $\pobs
(Y|X, \obsparam)$, could easily violate the structure of our model in
this case.
For instance, the same chemical parameters could affect both the
measurement and reaction processes, violating the assumed separation of
$\param$ and $\obsparam$.

Even careful preprocessing in such a setting can create a fundamental
incoherence.
Suppose the downstream analysis will be Bayesian, so the preprocessor
provides the conditional density of $Y$ as a function of $X$, $\pobs(Y
| X)$, for the observed $Y$.
If $\param$ and $\obsparam$ share components, and the preprocessor uses
their prior on $\obsparam$ to create $\pobs(Y | X)$, the conditional
density need not be sufficient for $\param$ under the downstream
analyst's model.
Because the downstream analyst's prior on $\param$ need not be
compatible with the preprocessor's prior on~$\obsparam$, inferences
based on the preprocessor's $\pobs(Y | X)$ can be seriously flawed in
this setting.
Hence, we exclude such cases from our investigation for the time being.

Thinking Bayesianly, our model (\ref{e:model}) obviously does not
exclude the possibility that the downstream analyst has more knowledge
about $\param$ than the preprocessor in the form of a prior on $\theta$.
However, \textit{prior} information means that it is based on studies
that do not overlap with the current one.
Probabilistically speaking, this means that our model permits the
downstream analyst to formally incorporate another data set $Z$, as
long as $Z$ is conditionally independent of the scientific variables
$X$ and observations $Y$ given $(\param, \obsparam)$ or~$\param$.
For example, the downstream analyst could observe completely separate
experiments pertaining to the same underlying process governed by
$\param$ or the outcomes of separate calibration pertaining to
$\obsparam$, but not additional replicates governed by the same
realization of $X$.
In a biological setting, this means that the downstream analyst could
have access to results from samples not available to the preprocessor
(e.g., biological replicates), possibly using the same equipment;
however, they could not have access to additional analyses of the same
biological sample (e.g., technical replicates), as a single biological
sample would typically correspond to a single realization of $X$.

These examples remind us that our multiphase setting does not encompass
all of statistical inference.
This is quite a relief to us.
Our work aims to open new directions for statistical research, but it
cannot possibly address every problem under the sun!

\subsection{Constraints will set your theory free}
\label{sec:constraints}

Multiphase theory hinges on procedural constraints.
Consider, for example, finding the optimal multiphase estimation
procedure in terms of the final estimator's Bayes risk.
Without stringent procedural constraints, the result is trivial:
compute the appropriate Bayes estimator using the distribution of $T$
given $\param$.
Similarly, the optimal preprocessing $T$ will, without tight
constraints, simply compute an optimal estimator using $Y$ and pass it forward.
Note that both of these cases respect risk-monotonicity to the letter;
it is not sufficiently tight to enable interesting, relevant theory.
More constraints, based upon careful consideration of applied problems,
are clearly required.

This is not altogether bad news.
We need only look to the history of multiple imputation to see how rich
theory can arise from stringent, pragmatic constraints.
Multiple imputation forms a narrow subset of multiphase procedures: $X$
corresponds to the complete data ($Y_{com}$, in MI notation), $Y$
corresponds to the observed data $Y_{obs}$ and missing data indicator
$R$, and $T$ usually consists of posterior predictive draws of the
missing data together with the observed data.
The Markovian property depicted in Figure~\ref{fig:models} holds when
the parameter ($\xi$) for the missing data mechanism $p(R|Y_{com},\xi$)
is distinct from the parameter of interest ($\theta$) in
$p(Y_{com}|\theta)$, which is a common assumption in practice.
The second-phase procedure is then restricted to repeatedly applying a
complete-data procedure and combining the results.
These constraints were originally imposed for practical reasons---in
particular, to make the resulting procedure feasible with existing software.
However, they have opened the door to deep theoretical investigations.

In that spirit, we consider two types of practically-motivated
constraints for multiphase inference: restrictions on the downstream
analyst's procedure and restrictions on the preprocessor's methods.
These constraints are intended to work in concert with coherence
conditions (e.g., risk monotonicity), not in isolation, to enable
meaningful theory.

Constraints on the downstream analyst are intended to reflect practical
limitations of their analytic capacity.
Examples include restricting the downstream analyst to narrow classes
of estimators (e.g., linear functions of preprocessed inputs), to
specific principles of estimation (e.g., MLEs), or to special cases of
a method we can reasonably assume the downstream analyst could handle,
such as a complete-data estimator $\est(X)$, available from software
with appropriate inputs.
Estimators derived from nested families of models are often suitable
for this purpose.
For example, whereas $\est(X)$ may involve only an ordinary regression,
the computation of $\est(T)$ may require a weighted least-squares
regression.\looseness=1

Another constraint on the downstream analyst pertains to nuisance parameters.
Such constraints are of great practical and theoretical interest, as we
believe that the preprocessor will typically have better knowledge and
statistical resources available to address nuisance parameters than the
downstream analyst.
An extreme but realistic case of this is to assume that the downstream
analyst cannot address nuisance parameters at all.
As we shall discuss in Section~\ref{sec:theory}, this would force the
preprocessor to either marginalize over the nuisance parameters, find a
pivot with respect to them, or trust the downstream analyst to use a
method robust to the problematic parameters.

Turning to the preprocessor, we consider restricting either the form of
the preprocessor's output or the mechanics of their methods.
In the simplest case of the former, we could require that $T$ consist
of the posterior mean ($\hat{X}$) and posterior covariance ($V$) of the
unknown $X$ under the preprocessor's model.
A richer, but still realistic, class of output would be
finite-dimensional real or integer vectors.
Restricting output to such a class would prevent the preprocessor from
passing arbitrary functions onto the downstream analyst.
This leads naturally to the investigation of (finite-dimensional)
approximations to the preprocessor's conditional density, aggregation,
and other such techniques.

On the mechanical side, we can restrict either the particulars of the
preprocessor's methods or their broader properties.
Examples of the former include particular computational approximations
to the likelihood function or restrictions to particular principles of
inference (e.g., summaries of the likelihood or posterior distribution
of $X$ given $(Y, \obsparam)$).
Such can focus our inquiries to specific, feasible methods of interest
or reflect the core statistical principles we believe the preprocessor
should take into account.
In a different vein, we can require that preprocessor's procedures be
distributable across multiple researchers, each with their own
experiments and scientific variables of interest.
Such settings are of interest for both the accumulation of scientific
results for later use and for the development of distributed
statistical computation.
This leads to preprocessing based upon factored ``working'' models for
$X$, as we explore further in Section~\ref{sec:sufficiency}.
Nuisance parameters play an important role in these constraints,
narrowing the class of feasible methods (e.g., marginalization over
such parameters may be exceedingly difficult) and largely determining
the extent to which preprocessing can be distributed.
We explore these issues in more detail throughout Section~\ref{sec:theory}.

\section{A few theoretical cornerstones}
\label{sec:theory}

We now present a few steps towards a theory of multiphase inference.
In this, we endeavor to address three basic questions: (1) how can we
determine what to retain, (2) what limits the performance of multiphase
procedures, and (3) what are some minimal requirements for being an
ideal preprocessor?
We find insight into the first question from the language of classical
sufficiency.
We leverage and specialize results from the missing-data literature to
address the second.
For the third question, we turn to the tools of decision theory.

\subsection{Determining what to retain}
\label{sec:sufficiency}

Suppose we have a group of researchers, each with their own experiments.
They want to preprocess their data to reduce storage requirements, ease
subsequent analyses, and (potentially) provide robustness to
measurement errors.
This group is keenly aware of the perils of preprocessing and want to
ensure that the output they provide will be maximally useful for later analyses.
Their question is, ``Which statistics should we retain?''

If each of these researchers was conducting the final analysis
themselves, using only their own data, they would be in a single-phase setting.
The optimal strategy then is to keep a minimal sufficient statistic for
each researcher's model.
Similarly, if the final analysis were planned and agreed upon among all
researchers, we would again have a single-phase setting, and it is
optimal to retain the sufficient statistics for the agreed-upon model.
We use the term \textit{optimal} here because it achieves maximal data
reduction without losing information about the parameters of interest.
Such lossless compression--in the general sense of avoiding statistical
redundancy--is often impractical, but it provides a useful theoretical
gold standard.

In the multiphase setting, especially with multiple researchers in the
first phase, achieving optimal preprocessing is far more complicated
even in theory.
If $T(Y)$ is the output of the \textit{entire} preprocessing phase,
then in order to retain all information we must require $T(Y)$ to be a
sufficient statistics for $\{\theta,\xi\}$ under model (\ref{e:model});
that is,
%
%
\begin{equation}
\label{e:con1} L\bigl(\theta,\xi|T(Y)\bigr)=L(\theta,\xi|Y),
\end{equation}
where $L$ denotes a likelihood function; or at least in the (marginal)
Bayesian sense,
%
%
\begin{equation}
\label{e:con2} P\bigl(\theta|T(Y)\bigr)=P(\theta|Y),
\end{equation}
where $P(\theta|D)$ is the posterior of $\theta$ given data $D$ with
the likelihood given by (\ref{e:model}).
Note that (\ref{e:con1}) implies (\ref{e:con2}), and (\ref{e:con2}) is
useful when the downstream analyst wants only a Bayesian inference of
$\theta$.
In either case the construction of the sufficient statistic generally
depends on the joint model for $Y$ as implied by (\ref{e:model}),
requiring more knowledge than individual researchers typically possess.

Often, however, it is reasonable to assume the following conditional
independence.
Let $\{Y_i, X_i, \obsparam_i\}$ be the specification of $\{Y, X,
\obsparam\}$ for researcher $i (=1,\ldots, r)$, where $\{Y_1,\ldots,
Y_r\}$ forms a \textit{partition} of $Y$. We then assume that
%
%
\begin{equation}
\label{e:obsm} \pobs(Y | X, \obsparam)= \prod_{i=1}^r
p_{Y_i}(Y_i | X_i, \obsparam _i) .
\end{equation}
Note in the above definition implicitly we also assume the baseline
measure $\mu_Y$ is a product measure $\prod_{i=1}^r \mu_{Y_i}$, such as
Lebesgue measure.
The assumption (\ref{e:obsm}) holds, for example, in microarray
applications, when different labs provide conditionally-independent
observations of probe-level intensities.
The preceding discussion suggests that this assumption is necessary for
ensuring (\ref{e:con1}) or even (\ref{e:con2}), but obviously it is far
from sufficient because it says nothing about the model on $X$.

It is reasonable---or at least more logical than not---to assume each
researcher has the best knowledge to specify his/her own observation
model $p_{Y_i}(Y_i | X_i, \obsparam_i)$ ($i=1,\ldots, r)$.
But, for the scientific model $p_X(X|\theta)$ used by the downstream
analyst, the best we can hope is that each researcher has \textit{a
working model} $\tilde p_X(X_{i}|g_i(\eta))$ that is in some way
related to $p_X(X|\theta)$.
The notation $g_i(\eta)$ reflects our hope to construct a common
working parameter $\eta$ that can ultimately be \textit{linked} to the
scientific parameter $\theta$.

Given this working model, the $i$th researcher can obtain the
corresponding (minimal) sufficient statistic $T_i(Y)$ for $\{
g_i(\workparam), \xi_i\}$ with respect to
%
%
\begin{equation}
\label{e:prob} \pwork\bigl(Y_i | g_i(\workparam),
\obsparam_i\bigr) = \int\pobs(Y_i | X_i,
\obsparam_i) \pwork\bigl(X_i | g_i(
\workparam)\bigr) \,\dd\mu_{X_i}(X_i),\qquad i=1,\ldots, r.
\end{equation}
When one has a prior $ \pi_{\obsparam_i}(\obsparam_i)$ for $\xi_i$, one
could alternately decide to retain the (Bayesian) sufficient statistic
$T_i^{B}(Y_{i})$ with respect to the model
%
%
\begin{equation}
\label{e:prom} \pywork\bigl(Y_i | g_i(\workparam)\bigr)
= \int\int\pobs(Y_i | X_i, \obsparam _i)
\pwork\bigl(X_i | g_i(\workparam)\bigr)
\pi_{\obsparam_i}(\obsparam_i) \,\dd\mu _{X_i}(X_i)
\,\dd\mu_\obsparam(\obsparam_i) .
\end{equation}

Our central interest here is to determine when the collection $T(Y)=\{
T_i(Y_i)\dvt i=1,\ldots, r\}$ will satisfy (\ref{e:con1}) and when
$T^{B}(Y)=\{T_i^B(Y_i)\dvt i=1,\ldots, r\}$ will satisfy (\ref{e:con2}).
This turns out to be an exceedingly difficult problem if we seek a
necessary and sufficient condition for \textit{when} this occurs.
However, it is not difficult to identify sufficient conditions that can
provide useful practical guidelines.
We proceed by first considering cases where $\{X_1, \ldots, X_r\}$
forms a partition of $X$.
Compared to the assumption on partitioning $Y$, this assumption is less
likely to hold in practice because different researchers can share
common parts of $X$'s or even the entire scientific variable $X$.
However, as we shall demonstrate shortly, we can extend our results
formally to all models for $X$, as long as we are willing to put tight
restrictions on the allowed class of working models.
Specifically, the following condition describes a class of working
models that are ideal because they permit separate preprocessing yet
retain joint information.
Note again that an implicit assumption here is that the baseline
measure $\mu_X$ is a product measure $\prod_{i=1}^r\mu_{X_i}$.
%
\begin{definition}
[(Distributed separability condition (DSC))]
A set of working models $\{\pwork(X_i | g_{i}(\workparam))\dvt i=1,\ldots,
r\}$ is said to satisfy the \emph{distributed separability condition}
with respect to $\psci(X | \param)$ if there exists a probability
measure $p_\workparam( \workparam| \param)$ such that
%
%
\begin{equation}
\psci(X | \param) = \int_\eta \Biggl[ \prod
_{i=1}^r \pwork\bigl(X_i |
g_i(\workparam)\bigr) \Biggr]\,\dd p_\workparam(\workparam|
\param). \label{eq:dsc}
\end{equation}
\end{definition}


\begin{theorem}
\label{thm:dsc}
Under the assumptions (\ref{e:obsm}) and (\ref{eq:dsc}), we have
\begin{longlist}[(1)]
\item[(1)] The collection of individual sufficient statistics from
(\ref
{e:prob}), that is, $T(Y)=\{T_i(Y_i), i=1,\ldots, r\}$, is jointly
sufficient for $\{\theta,\xi\}$ in the sense that (\ref{e:con1}) holds.
\item[(2)] Under the additional assumption that $\{\xi_1,\ldots, \xi
_r\}
$ forms a partition of $\xi$ and $\pi(\xi)\,\dd\mu_{\xi}=\prod_{i=1}^r
\pi_{\xi_i}(\xi_i)\,\dd\mu_{\xi_i}$, both $T(Y)$ corresponding to
(\ref
{e:prob}) and $T^B(Y)$ corresponding to (\ref{e:prom}) are Bayesianly
sufficient for $\theta$ in the sense that (\ref{e:con2}) holds.
\end{longlist}
\end{theorem}
%
%
\begin{pf}
By the sufficiency of $T_i$ for $(g_i(\workparam), \obsparam_i)$, we
can write
%
%
\begin{equation}
\label{eq:a1} \int_{X_i} \pobs(Y_i |
X_i, \obsparam_i) \pwork\bigl(X_i |
g_i(\workparam)\bigr) \,\dd\mu_{X_i}(X_i) =
\pywork\bigl(Y_i | g_i(\workparam),
\obsparam_i\bigr) =h_i(Y_i)
f_i\bigl(T_i ; g_i(\workparam),
\obsparam_i\bigr).
\end{equation}
This implies that,
\begin{eqnarray*}
\pobs(Y | \param, \obsparam) &=& \int_{X} \pobs(Y | X,
\obsparam) \psci (X | \param)\,\dd\mu_X(X),
\\
{\bigl[\mbox{by } (\ref{e:obsm})  \mbox{ and } (\ref{eq:dsc})\bigr]}& =& \int
_{X} \Biggl[ \prod_{i=1}^{r}
\pobs(Y_i | X_i, \obsparam_i) \Biggr]
\\
&&{} \times\Biggl[ \int_{\eta} \Biggl[ \prod
_{i=1}^{r} \pwork\bigl(X_i |
g_i(\workparam )\bigr) \Biggr]\,\dd p_\workparam(\workparam|
\param) \Biggr]\,\dd\mu _X(X),
\\
{[\mbox{by factorization of }  \mu_X]} &=& \int
_{\eta} \prod_{i=1}^{r}
\biggl[ \int_{X_i} \pobs(Y_i | X_i,
\obsparam_i) \pwork\bigl(X_i | g_i(
\workparam)\bigr)\,\dd\mu_{X_i}(X_i) \biggr]\,\dd
p_\workparam (\workparam| \param),
\\
{\bigl[\mbox{by } (\ref{eq:a1})\bigr]} &=& \Biggl[ \prod
_{i=1}^{r} h_i(Y_i) \Biggr]
\Biggl[ \int_{\eta} \prod_{i=1}^{r}
f_i\bigl(T_i;g_i(\workparam), \obsparam
_i\bigr)\,\dd p_\workparam(\workparam| \param) \Biggr] .
\end{eqnarray*}
This establishes (1) by the factorization theorem.
Assertion (2) is easily established via an analogous argument, by
integrating all the expressions above with respect to $\pi(\xi)\,\dd
\mu
(\xi)=\prod_{i=1}^r \pi_{\xi_i}(\xi_i)\,\dd\mu_{\xi_i}(\xi_i)$.
\end{pf}

We emphasize that DSC does not require individual researchers to model
their parts of $X$ in the same way as the downstream analyst would,
which would make it an essentially tautological condition.
Rather, it requires that individual researchers understand their own
problems and how they can fit into the broader analysis hierarchically.
This means that the working model for each $X_i$ $(i=1,\ldots, r)$ can
be more saturated than the downstream analyst's model for the same part
of $X$.

Consider a simple case with $r=1$, where the preprocessor correctly
assumes the multivariate normality for $X$ but is unaware that its
covariance actually has a block structure or is unwilling to impose
such a restriction to allow for more flexible downstream analyses.
Clearly any sufficient statistic under the unstructured multivariate
model is also sufficient for any (nested) structured ones.
The price paid here is failing to achieve the greatest possible
sufficient reduction of the data, but this sacrifice may be necessary
to ensure the broader validity of downstream analyses.
For example, even if downstream analysts adopt a block-structured
covariance, they may still want to perform a model checking, which
would not be possible if all they are given is a \textit{minimal}
sufficient statistic for the model to be checked.

Knowledge suitable for specifying a saturated model is more attainable
than complete knowledge of $\psci(X | \param)$, although ensuring
common knowledge of its (potential) hierarchical structure still
requires some coordination among the researchers.
Each of them could independently determine for which classes of
scientific models their working model satisfies the DSC.
However, without knowledge of the partition of $X$ across researchers
and the overarching model(s) of interest, their evaluations need not
provide any useful consensus.
This suggests the necessity of some general communications and a
practical guideline for distributed preprocessing, even when we have
chosen a wise division of labors that permits DSC to hold.

Formally, DSC is similar in flavor to de Finetti's theorem, but it does
not require the components of the factorized working model to be exchangeable.
DSC, however, is by no means necessary (even under (\ref{e:con1})), as
an example in Section~\ref{sec:counterexamples} will demonstrate.
Its limits stem from ``unparameterized'' dependence---dependence
between $X_i$'s that is not controlled by $\param$.
When such dependence is present, statistics can exist that are
sufficient for both $\workparam$ and $\param$ without the working model
satisfying DSC.

However, a simple necessary condition for distributed sufficiency is available.
Unsurprisingly, it links the joint sufficiency of $T(Y)=\{T_i(Y_i):
i=1,\ldots, r\}$ under $p_{Y}(Y|\theta)$ to the joint sufficiency of
$S(X)=\{S_i(X_i), i=1,\ldots,r\}$ under the scientific model
$p_X(X|\theta)$, where $S_i(X_i)$ is any sufficient statistic for the
working model $\tilde p_X(X_{i}|g_i(\eta)), i=1,\ldots, r$.
%
\begin{theorem}
\label{thm:necessary}
If, for all observation models satisfying (\ref{e:obsm}), the
collections of individual sufficient statistics from (\ref{e:prob})
$T(Y)=\{T_i(Y_i), i=1,\ldots, r\}$ are jointly sufficient for $\{
\theta
,\xi\}$ in the sense that (\ref{e:con1}) holds, then any collection of
individual sufficient statistics under $\{\pwork(X_i | g_i(\eta)),
i=1,\ldots, r\}$, that is, $S(X)$, must be sufficient for $\theta$
under $\psci(X | \theta)$.
\end{theorem}
The proof of this condition emerges easily by considering the trivial
observation model $\pobs(Y_i | X_i, \xi) = \delta_{\{Y_i = X_i\}}$,
where $\delta_{A}$ is the indicator function of set $A$.
Theorem~\ref{thm:necessary} holds even if we require the observation
model to be nontrivial, as the case of $\pobs(Y_i | X_i, \xi) \propto
\delta_{\{Y_i \in B_{\varepsilon}[X_i]\}}$ for arbitrary $\varepsilon
$-neighborhoods of $X_i$ demonstrates.
The result says that if we want distributed preprocessing to provide a
lossless compression regardless of the actual form of the observation
model, then even under the conditional independence assumption (\ref
{e:obsm}), we must require the individual working models to \textit
{collectively} preserve sufficiency under the scientific model.
Note that preserving sufficiency for a model is a much weaker
requirement than preserving the model itself.
Indeed, two models can have very different model spaces yet share the
same \textit{form} of sufficient statistics, as seen with i.i.d.
$\Poisson(\mu)$ and $N(\mu, 1)$ models, both yielding the sample
average as a complete sufficient statistic.

Although we find this sufficiency-preserving condition quite
informative about the limits of lossless distributed preprocessing, it
is not a sufficient condition.
As a counterexample, consider $Y_{ij} | X_i \sim N(\mu_i, \sigma_{ij}^2
)$ independent for $i=1,\ldots,n$, $j=1,\ldots,m$, where $X_i \equiv
(\mu_i, \sigma_{i1}^2, \ldots, \sigma_{im}^2)$.
For the true model, we assume $\psci(X|\theta)$ as
follows: $\mu_i | \theta\sim N(\theta, 1)$, $\sigma_{ij}^2 \sim
1/{\chi
_1}^2$, and all variables are mutually independent.
For the working model, we take $\pwork(X|\eta)$ as
follows: $\mu_i | \eta_i \sim N(\eta_i, 1)$ independently, and
$\sigma
_{ij}^2 = 1$ with probability 1 for all $i,j$.
Obviously $S=(\mu_1, \ldots, \mu_n)$ is a sufficient statistic for both
$\pwork(X|\eta)$ and $p_X(X|\theta)$ because of their normality.
Because $S$ is \textit{minimally} sufficient for $\eta$, this implies
that any sufficient statistic for $\pwork(X|\eta)$ must be sufficient
for $p_X(X|\theta)$, therefore the sufficiency preserving condition holds.

However, the collection of the complete sufficient statistics $T_i =
\sum_j y_{ij}/m, i=1, \ldots, r$ for $\eta$ under $p_Y(Y|\eta)$ is not
sufficient for $\theta$ under $\pobs(Y | \theta)$ because the latter is
no longer an exponential family.
The trouble is caused by the failure of the working models to capture
additional flexibility in the scientific model that is not controlled
by its parameter $\theta$.
Therefore, obtaining a condition that is both necessary and sufficient
for lossless compression via distributed preprocessing is a challenging task.
Such a condition appears substantially more intricate than those
presented in Theorems~\ref{thm:dsc} and~\ref{thm:necessary} and may
therefore be less useful as an applied guideline.
Below we discuss a few further subtleties.

\subsubsection*{Likelihood sufficiency verses Bayesian efficiency}
Although Theorem~\ref{thm:dsc} covers both likelihood and Bayesian
cases, it is important to note a subtle distinction between their
general implications. In the likelihood setting (\ref{e:con1}), we
achieve lossless compression for all downstream analyses targeting
$(\param, \obsparam)$. This allows the downstream analyst to obtain
inferences that are robust to the preprocessor's beliefs about
$\obsparam$, and they are free to revise their inferences if new
information about $\obsparam$ becomes available. But, the downstream
analyst must address the nuisance parameter $\xi$ from the
preprocessing step, a task a downstream analyst may not be able or
willing to handle.

In contrast, the downstream analyst need not worry about $\xi$ in the
Bayesian setting (\ref{e:con2}).
However, this is achieved at the cost of robustness.
All downstream analyses are potentially affected by the preprocessors'
beliefs about $\xi$.
Furthermore, because $T^B(Y)$ is required only to be sufficient for
$\theta$, it may not carry any information for a downstream analyst to
check the preprocessor's assumptions about $\obsparam$.
Fortunately, as it is generally logical to expect the preprocessor to
have better knowledge addressing $\xi$ than the downstream analyst,
such robustness may not be a serious concern from a practical perspective.
Theoretically, the trade-off between robustness and convenience is not
clear-cut; they can coincide for other types of preprocessing, as seen
in Section~\ref{sec:missinfo} below.

\subsubsection*{Deterministic dependencies among $X_i$'s}
As discussed earlier, (conditional) dependencies among the observation
variables $Y_i$ across different $i$'s will generally rule out the
possibility of achieving lossless compression by collecting individual
sufficient statistics. This points to the importance of appropriate
separation of labors when designing distributed preprocessing.
In contrast, dependencies among $X_i$'s are permitted, at the expense
of redundancy in sufficient statistics.
We first consider deterministic dependencies, and for simplicity, take
$r = 2$ and constrain attention to the case of sufficiency for $\param$.
Suppose we have $X_1$ and $X_2$ forming a partition of $X$, with a
working model $\pwork(X | \bm\workparam)=\tilde p_{X_1}(X_1 |\eta
_1)\tilde p_{X_2}(X_2 |\eta_2) $ that satisfied the DSC for some
$p_\workparam(\workparam| \param)$.
Imagine we need to add a common variable $Z$ to both $X_1$ and $X_2$
that is conditionally independent of $\{X_1, X_2\}$ given $\theta$ and
has density $p_Z(Z|\theta)$, with the remaining model unchanged.
However, the two researchers are unaware of the sharing of~$Z$, so they
set up $X_1'=\{X_1, Z_1\}$ and $X_2'=\{X_2, Z_2\}$, with $p_{X_1'}(X_1'
|\eta_1')$ and $\tilde p_{X_2'}(X_2' |\eta_2')$ as their respective
working models.\vspace*{1pt}

At the first sight this seems to be a hopeless situation for applying
the DSC condition, because $X'=\{X_1', X_2'\}=\{X_1, Z_1, X_2, Z_2\}$
does not correspond to the scientific variable $X=\{X_1, X_2, Z\}$ of interest.
However, we notice that if we can force $Z_1=Z_2=Z$ in $X'$, then we
can recover $X$.
This forcing is not a mere mathematical trick.
Rather, it reflects an extreme yet practical strategy when researchers
are unsure whether they share some components of their $X_i's$ with others.
The strategy is simply to retain statistics sufficient for the entire
part that they may \textit{suspect} to be common, which in this case
means that both researchers will retain statistics sufficient for the
$Z_i'$s $(i=1,2)$ in their entirety.
Mathematically, this corresponds to letting $\tilde p_{X_i'}(X_i' |\eta
_i')= \tilde p_{X_i}(X_i |\eta_i)\delta_{\{Z_i=\zeta_i\}}$, where
$\eta
'_i=\{\eta_i,\zeta_i\}$.
It is then easy to verify that DSC holds, if we take
$p_\workparam'(\workparam' | \param)=p_\workparam(\workparam|
\param
)p_{Z}(\zeta_1|\theta)\delta_{\{\zeta_1=\zeta_2\}}$, where $\eta
'=\{\eta
, \zeta_1, \zeta_2\}$. This is because when $Z_1\neq Z_2$, both sides
of (\ref{eq:dsc}) are zero. When $Z_1=Z_2=Z$, we have (adopting
integration over $\delta$ functions)
\begin{eqnarray*}
&&\int_{\eta'} \Biggl[ \prod_{i=1}^2
\tilde p_{X_i'}\bigl(X_i' |
\eta_i'\bigr) \Biggr]\,\dd p_\workparam'
\bigl( \workparam' | \param\bigr) \\
&&\quad= \int_\eta
\int_{\zeta
_1} \Biggl[ \prod_{i=1}^2
\tilde p_{X_i}(X_i |\eta_i)
\delta_{\{
Z=\zeta_i\}
} \Biggr]\,\dd p_\workparam( \workparam| \param)
\delta_{\{\zeta_1=\zeta_2\}}\,\dd p_{Z}(\zeta_1|\theta)
\\
&&\quad = \Biggl[ \int_\eta\prod_{i=1}^2
\tilde p_{X_i}(X_i |\eta_i)\,\dd
p_\workparam( \workparam| \param) \Biggr] \int_{\zeta_1}
\delta_{\{\zeta_1=Z\}}\,\dd p_{Z}(\zeta_1|\theta)
\\
&&\quad =p_{X_{}}(X_1, X_2|\theta)p_{Z}(Z|
\theta)=p_{X}(X|\theta).
\end{eqnarray*}
This technique of expanding $\eta$ to include shared parts of the $X$
allows the DSC and Theorem~\ref{thm:dsc} to be applied to all models
$\psci(X | \theta)$, not only those with with distinct $X_i$'s.
However, this construction also restricts working models to those with
deterministic relationships between parts of $\eta$ and each $X_i$.

The derivation above demonstrates both the broader applications of DSC
as a theoretical condition and its restrictive nature as a practical
guideline. Retaining sufficient statistics for both $Z_1$ and $Z_2$ can
create redundancy.
If each preprocessor observes $Z$ without noise, then only one of them
actually needs to retain and report their observation of $Z$.
However, if each observes $Z$ with independent noise, then both of
their observations are required to obtain a sufficient statistic for~$\theta$.
The noise-free case also provides a straightforward counterexample to
the necessity of DSC. Assuming both preprocessors observe $Z$ directly,
as long as one of the copies of $Z$ is retained via the use of the
saturated $\delta$ density, the other copy can be modeled in any way---and
hence can be made to violate DSC---without affecting their joint
sufficiency for~$\theta$.

Regardless of the dependencies among the $X_i$'s, there is always a
safe option open to the preprocessors for data reduction: retain $T_i$
sufficient for $(X_i, \obsparam_i)$ under $\pobs(Y_i | X_i, \obsparam_i)$.
This will preserve sufficiency for $\param$ under any scientific model
$\psci(X | \param)$:
%
\begin{theorem}\label{thm:safe}
If $\pobs(Y | X, \obsparam)$ is correctly specified and satisfies
(\ref
{e:obsm}), then any collection of individual sufficient statistics $\{
T_i : i=1,\ldots,r\}$ with each $T_i$ sufficient for $(X_i, \obsparam
_i)$ is jointly sufficient for $(\param, \obsparam)$ in the sense of
(\ref{e:con1}) for all $\psci(X | \param)$.
\end{theorem}
\begin{pf}
By the factorization theorem, we have $\pobs(Y_i | X_i, \obsparam_i) =
h_i(Y_i) f_i(T_i ; X_i, \obsparam_i)$ for any $i$. Hence, by (\ref
{e:obsm}), $ \pobs(Y | \theta) =  [ \prod_{i=1}^r h_i(Y_i)
]
\int_X  [ \prod_{i=1}^r p_T(T_i | X_i, \obsparam_i)  ]
\psci(X
| \theta) \,\dd\mu_X( X)$. Therefore $\{T_i : i=1 ,\ldots,r\}$ is
sufficient for $\theta$, by the factorization theorem for sufficiency.
\end{pf}
Theorem~\ref{thm:safe} provides a universal, safe strategy for
sufficient preprocessing and a lower bound on the compression
attainable from distributed sufficient preprocessing.
As all minimal sufficient statistics for $\theta$ are functions of any
sufficient statistic for $(X, \obsparam)$, retaining minimal sufficient
statistics for each $(X_i, \obsparam_i)$ results in less compression
than any approach properly using knowledge of $\psci(X | \theta)$.
However, the compression achieved relative to retaining $Y$ itself may
still be significant.
Minimal sufficient statistics for $\theta$ provide an upper bound on
the attainable degree of compression by the same argument.
Achieving this compression generally requires that each preprocessor
knows the true scientific model $\psci(X | \theta)$.
Between these bounds, the DSC (\ref{eq:dsc}) shows a trade-off between
the generality of preprocessing (with respect to different scientific
models) and the compression achieved: the smaller the set of scientific
models for which a given working model satisfies (\ref{eq:dsc}), the
greater the potential compression from its sufficient statistics.

\subsubsection*{Stochastic dependencies among $X_i$'s}
More generally, stochastic
dependence among $X_i$'s reduces compression and increases redundancy
in distributed preprocessing.
These costs are particularly acute when elements of $\param$ control
dependence among $X_i$'s, as seen in the following example where
\begin{eqnarray*}
X &=& (X_1, X_2)^\top \sim
N_{4D} \left(\param_1 \bm1_{4D}, \pmatrix{ 1 & 0 & 0
& \param_2 \vspace*{2pt}
\cr
0 & 1 & 0 & 0 \vspace*{2pt}
\cr
0 & 0 & 1
& 0 \vspace*{2pt}
\cr
\param_2 & 0 & 0 & 1 } \otimes I_{D}
\right)\qquad\mbox{for }D > 1,
\\
Y_i &= &(Y_{i1}, Y_{i2})^\top
\given X_i\sim N_{2D}(X_i, I_{2D})
\qquad\mbox{independently for } i=1,2.
\end{eqnarray*}
Here $\bm1_{4D}$ is a column vector with $4D$ $1$'s as its components,
and $\otimes$ is the usual Kronecker product.
If $\param_2$ is known, then each researcher can reduce their
observations $Y_i$ to a scalar statistic $Y_i^\top\bm1_{2D}$ and
preserve sufficiency for $\param_1$.
If $\param_2$ is unknown, then each researcher must retain all of
$Y_{ii}$ (but not $Y_{ij}$ for $i\not=j$) in addition to these sums to
ensure sufficiency for $\param= (\param_1, \param_2)$, because the
minimal sufficient statistic for $(\param_1, \param_2)$ requires the
computation of $Y_{11}^\top Y_{22}$.
Thus, the cost of dependence here is $D$ additional pieces of
information per preprocessor.
Dependence among the $X_i$'s forces the preprocessors to retain enough
information to properly combine their individual contributions in the
final analysis, downweighting redundant information.
This is true even if they are interested only in efficient estimation
of $\param_1$, leading to less reduction of their raw data and less
compression from preprocessing than the independent case.

\subsubsection*{Practical perspective}
From this investigation, we see that it is generally not enough for
each researcher involved in preprocessing to reduce data based on even
a correctly-specified model for their problem at hand.
We instead need to look to other models that include each
experimenter's data hierarchically, explicitly considering higher-level
structure and relationships.
However, significant reductions of the data are still possible despite
these limitations.
Each $T_i$ need not be sufficient for each~$X_i$, nor must $T$ be
sufficient for $X$ overall.
This often implies that much less data need to be retained and shared
than retaining sufficient statistics for each $X_i$ would demand.
For instance, if a working model with $X_i | \workparam_i \sim N(\mu_i,
\Sigma_i)$ satisfies the DSC for a given model $\psci(X | \param)$ and
$Y_{ij} | X_i, \obsparam_i \sim N(X_i, \obsparam_i)$, then only means
and covariance matrices of $Y_{ij}$ within each experiment $i$ need to
be retained.

The discussions above demonstrate the importance of involving
downstream analysts in the design of preprocessing techniques.
Their knowledge of $\psci(X | \param)$ is extremely useful in
determining what compression is appropriate, even if said knowledge is
imperfect.
Constraining the scientific model to a broad class may be enough to
guarantee effective preprocessing.
For example, suppose we fix a working model and consider all scientific
models that can be expressed as~(\ref{eq:dsc}) by varying the choices
of $p_\workparam( \workparam| \param)$.
This yields a very broad class of hierarchical scientific models for
downstream analysts to evaluate, while permitting effective distributed
preprocessing based on the given working model.\looseness=1

Practically, we see two paths to distributed preprocessing:
coordination and caution.
Coordination refers to the downstream analyst evaluating and guiding
the design of preprocessing as needed.
Such guidance can guarantee that preprocessed outputs will be as
compact and useful as possible.
However, it is not always feasible.
It may be possible to specify preprocessing in detail in some
industrial and purely computational settings.
Accomplishing the same in academic research or for any research
conducted over time is an impractical goal.
Without such overall coordination, caution is needed.
It is not generally possible to maintain sufficiency for $\param$
without knowledge of the possible models $\psci(X | \param)$ unless the
retained summaries are sufficient for $X$ itself.
Preprocessors should therefore proceed cautiously, carefully
considering which scientific models they effectively exclude through
their preprocessing choices.
This is analogous to the oft-repeated guidance to include as many
covariates and interactions as possible in imputation models (\citet
{Meng1994}, \citet{Meng2003}).

\subsection{Doing the best with what you get}
\label{sec:missinfo}

Having considered the lossless preprocessing, we now turn to more
realistic but less clear-cut situations.
We consider a less careful preprocessor and a sophisticated downstream analyst.
The preprocessor selects an output $T$, which may discard much
information in $Y$ but nevertheless preserves the identifiability of
$\theta$, and the downstream analyst knows enough to make the best of
whatever output they are given.
That is, the index set $C$ completely and accurately captures all
relevant preprocessing methods $T=\{ T_i\dvt i=1,\ldots, r \}$.
This does not completely capture all the practical constraints
discussed in Section~\ref{sec:concepts}.
However, it is important to establish an upper bound on the performance of multiphase procedures before incorporating such issues.
This upper bound is on the Fisher information, and hence a lower bound
on the asymptotic variances of estimators $\est$ of $\param$.
As we will see, nuisance parameters ($\obsparam$) play a crucial role
in these investigations.

When using a lossy compression, an obvious question is how much
information is lost compared to a lossless compression. This question
has a standard asymptotic answer when the downstream analyst adopts an
MLE or Bayes estimator, so long as nuisance parameters behave
appropriately (as will be discussed shortly). If the downstream analyst
adopts some other procedures, such as an estimating equation, then
there is no guarantee that the procedure based on $Y$ is more efficient
than the one based on $T$. That is, one can actually obtain a more
efficient estimator with less data when one is not using \textit
{probabilistically principled} methods, as discussed in detail in
\citet
{Meng2012}.

Therefore, as a first step in our theoretical investigations, we will
focus on MLEs; the results also apply to Bayesian estimators under the
usual regularity conditions to guarantee the asymptotic equivalence
between MLEs and Bayesian estimators. Specifically, let $(\hat{\param
}(Y), \hat{\obsparam}(Y))$ and $(\hat{\param}(T), \hat{\obsparam}(T))$
be the MLEs of $(\theta, \xi)$ based respectively on $Y$ and $T$ under
model (\ref{e:model}). We place standard regularity conditions for the
joint likelihood of $(\param, \obsparam)$, assuming bounded third
derivatives of the log-likelihood, common supports of the observation
distributions with respect to $(\param, \obsparam)$, full rank for all
information matrices at the true parameter value $(\theta_0, \xi_0)$,
and the existence of an open subset of the parameter space that
contains $(\theta_0, \xi_0)$. These conditions imply the first and
second Bartlett identities.

However, the most crucial assumption here is a sufficient accumulation
of information, indexed by an \textit{information size} $N_Y$, to
constrain the behavior of remainder terms in quadratic approximations
of the relevant score functions.
Independent identically distributed observations and fixed-dimensional
parameters would satisfy this requirement, in which case $N_Y$ is
simply the data size of $Y$, but weaker conditions can suffice (for an overview, see \citet{Lehmann1998}).
In general, this assumption requires that the dimension of both $\theta
$ and $\xi$ are bounded as we accumulate more data, preventing the type
of phenomenon revealed in \citet{NeymanScott1948}. For multiphase
inferences, cases where these dimensions are unbounded are common (at
least in theory) and represent interesting settings where preprocessing
can actually improve asymptotic efficiency, as we discuss shortly.

To eliminate the nuisance parameter $\xi$, we work with the observed
Fisher information matrices based on the profile likelihoods for
$\param
$, denoted by $I_Y$ and $I_T$ respectively. Let $F$ be the limit of
$I_Y^{-1}(I_Y-I_T)$, the so-called \textit{fraction of missing
information} (see \citet{Meng1991}), as $N_Y \rightarrow\infty$.
The proof of the following result follows the standard asymptotic
arguments for MLEs, with the small twist of applying them to profile
likelihoods instead of full likelihoods. (We can also invoke the more
general arguments based on decomposing estimating equations, as given
in \citet{Xie2012}.)
%
\begin{theorem}
\label{thm:missinfo} Under the conditions given above, we have
asymptotically as $N_Y \rightarrow\infty$,
%
%
\begin{equation}
\label{eq:regret} \Var \bigl(\est(T) - \est(Y) \bigr) \bigl[\Var \bigl(\est (T)
\bigr) \bigr]^{-1} \rightarrow F
\end{equation}
and
%
%
\begin{equation}
\label{eq:frac} \Var \bigl(\est(Y) \bigr) \bigl[\Var \bigl(\est(T) \bigr)
\bigr]^{-1} \rightarrow I - F.
\end{equation}
\end{theorem}

This establishes the central role of the fraction of missing
information $F$ in determining the asymptotic efficiency of multiphase
procedures under the usual asymptotic regime.
As mentioned above, this is an ideal-case bound on the performance of multiphase procedures, and it is based on the usual squared-error loss; both the
asymptotic regime and amount of knowledge held by the downstream
analyst are optimistic.
We explore these issues below, focusing on (1) mutual knowledge and
alternative definitions of efficiency, (2) the role of
reparameterization, (3) asymptotic regimes and multiphase efficiency,
and (4) the issue of robustness in multiphase inference.

\subsubsection*{Mutual knowledge and efficiency}
In practice, downstream analysts are unlikely to have complete
knowledge of $\pobs$. Therefore, even if they were given the entire
$Y$, they would not be able to produce the optimal estimator $\hat
\theta
(Y)$, making
the $F$ value given by Theorem~\ref{thm:missinfo} an unrealistic
yardstick. Nevertheless, Theorem~\ref{thm:missinfo} suggests a
direction for a more realistic standard.

The classical theory of estimation focuses on losses of the form
$L(\est
, \param_0)$, where $\theta_0$ denotes the truth.
Risk based on this type of loss, given by $R(\est, \param_0) =
E[L(\est
, \param_0)]$, is a raw measure of performance, using the truth as a baseline.
An alternative is regret, the difference between the risk of a given
estimator and an ideal estimator $\est^*$; that is, $R(\est, \param_0)
- R(\est^*, \param_0)$.
Regret is popular in the learning theory community and forms the basis
for oracle inequalities.
It provides a more adaptive baseline for comparison than raw risk, but
we can push further.
Consider evaluating loss with respect to an estimator rather than the truth.
For mean-squared error, this yields
%
%
\begin{equation}
\label{eq:risk} R \bigl(\est(T), \est(Y) \bigr)=E \bigl[ \bigl(\est(T) - \est (Y)
\bigr)^{\top} \bigl(\est(T) - \est(Y) \bigr) \bigr] .
\end{equation}
Can this provide a better baseline, and what are its properties?

For MLEs, $R(\est(T), \est(Y))$ behaves the same (asymptotically) as
additive regret because Theorem~\ref{thm:missinfo} implies that, as
$N_Y \rightarrow\infty$ under the classical asymptotic regime,
%
%
\begin{eqnarray}
\label{eq:same} R\bigl(\est(T), \est(Y)\bigr)&=&\Var\bigl(\est(T) - \est(Y)\bigr)=
\Var\bigl(\est(T)\bigr) -\Var \bigl(\est(Y)\bigr)
\nonumber
\\[-8pt]
\\[-8pt]
\nonumber
&=&R\bigl(\est(T),
\param_0\bigr)-R\bigl(\est(Y), \param_0\bigr) .
\end{eqnarray}
For inefficient estimators, (\ref{eq:same}) does not hold in general
because $\est(T) - \est(Y)$ is no longer guaranteed to be
asymptotically uncorrelated with $\est(Y)$.
In such cases, this is precisely the reason $\est(T)$ can be more
efficient than $\est(Y)$ or, more generally, there exists a constant
$\lambda\neq0$ such that $\lambda\est(T)+(1-\lambda)\est(Y)$ is
(asymptotically) more efficient than $\est(Y)$.
In the terminology of \citet{Meng1994}, the estimation procedure $\est
(\cdot)$ is not \textit{self-efficient} if (\ref{eq:same}) does not
hold, viewing $Y$ as the complete data $Y_{\mathrm{com}}$ and $T$ as the
observed data~$Y_{\mathrm{obs}}$.
Indeed, if $R(\est(T), \param_0) < R(\est(Y), \param_0)$, $R(\est(T),
\est(Y))$ may actually be \emph{larger} for a \emph{better} $\est(T)$
because of the inappropriate baseline $\est(Y)$; it is a measure of
difference, not dominance, in such cases.
Hence, some care is needed in interpreting this measure.

Therefore, we can view (\ref{eq:risk}) as a generalization of the usual
notion of regret, or the relative regret if we divide it by $R(\est(Y),
\param_0)$. This generalization is appealing for the study of
preprocessing: we are evaluating the estimator based on preprocessed
data directly against what could be done with the complete raw data,
sample by sample, and we no longer need to impose the restriction that
the downstream analysts must carry out the most efficient estimation
under a model that captures the actual preprocessing. This direction is
closely related to the idea of strong efficiency from \citet{Xie2012}
and \citet{Meng2012}, which generalizes the idea of asymptotic
decorrelation beyond the simple (but instructive) setting covered here.
Such ideas from the theory of missing data provide a strong
underpinning for the study of multiphase inference and preprocessing.

\subsubsection*{Reparameterization}
Theorem~\ref{thm:missinfo} also emphasizes the range of effects that
preprocessing can have, even in ideal cases.
Consider the role that $F$ plays under different transformations of
$\param$.
Although the eigenvalues of $F$ are invariant under one-to-one
transformations of the parameters, submatrices of $F$ can change substantially.
Formally, if $\param= (\param_1, \param_2)$ is transformed to
$\omega
= (\omega_1, \omega_2) = (g_1(\param_1, \param_2), g_2(\param_1,
\param
_2))$, then the fraction of missing information for $\omega_1$ can be
very different from that for $\param_1$.
These changes mean that changes in parameterization can reallocate the
fractions of missing information among resulting subparameters in
unexpected---and sometimes very unpleasant---ways.
This is true even for linear transformations; a given preprocessing
technique can preserve efficiency for $\param_1$ and $\param_2$
individually while performing poorly for $\param_1 - \param_2$.
Such issues have arisen in, for instance, the work of \citet{Xie2012}
when attempting to characterize the behavior of multiple imputation
estimators under uncongeniality.

\subsubsection*{Asymptotic regimes and multiphase efficiency}
On a fundamental level, Theorem~\ref{thm:missinfo} is a negative result
for preprocessing, at least for MLEs.
Reducing the data from $Y$ to $T$ can only hinder the downstream analyst.
Formally, this means that $I_T \leq I_Y$ (asymptotically) in the sense
that $I_Y - I_T$ is positive semi-definite.
As a result, $\est(Y)$ will dominate $\est(T)$ in asymptotic variance
for any preprocessing $T$.
Thus, the only justification for preprocessing appears to be pragmatic;
if the downstream analyst could not make use of $\pobs$ for efficient
inference or such knowledge could not be effectively transmitted,
preprocessing provides a feasible way to obtain the inferences of interest.
However, this conclusion depends crucially on the assumed behavior of
the nuisance parameter $\obsparam$.

The usual asymptotic regime is not realistic for many multiphase
settings, particularly with regards to $\obsparam$.
In many problems of interest, $\dim(\obsparam) / N_Y $ does not tend to
zero as $N_Y$ increases, preventing sufficient accumulation of
information on the nuisance parameter $\xi$.
A typical regime of this type would accumulate observations $Y_i$ from
individual experiments $i$, each of which brings its own nuisance
parameter $\obsparam_i$.
Such a process could describe the accumulation of data from
microarrays, for instance, with each experiment corresponding to a chip
with its own observation parameters, or the growth of astronomical
datasets with time-varying calibration.
In such a regime, preprocessing can have much more dramatic effects on
asymptotic efficiency.

In the presence of nuisance parameters, inference based on $T$ can be
more robust and even more efficient than inference based on $Y$.
It is well-known that the MLE can be inefficient and even inconsistent
in regimes where $\dim(\obsparam) \rightarrow\infty$ (going back
to at least \citet{NeymanScott1948}).
Bayesian methods provide no panacea either.
Marginalization over the nuisance parameter $\obsparam$ is appealing,
but resulting inferences are typically sensitive to the prior on
$\obsparam$, even asymptotically.
In many cases (such as the canonical Neyman--Scott problem), only a
minimal set of priors provide even consistent Bayes estimators.
Careful preprocessing can, however, enable principled inference in such regimes.

Such phenomena stand in stark contrast to the theory of multiple imputation.
In that theory, complete data inferences are typically assumed to be
valid. 
Thus, under traditional missing data mechanisms, the observed data
(corresponding to~$T$) cannot provide better inferences than $Y$.
This is not necessarily true in multiphase settings.
If the downstream analyst is constrained to particular principles of
inference (e.g., MLE or Bayes), then estimators based on $T$ can
provide lower asymptotic variance than those based on~$Y$.
This occurs, in part, because the mechanisms generating $Y$ and $T$
from $X$ are less restricted in the multiphase setting compared to the
traditional missing-data framework.
Principled inferences based on $X$ would, in the multiphase setting,
generally dominate those based on either $Y$ or $T$.
However, such a relationship need not hold between $Y$ and $T$ without
restrictions on the behavior of $\obsparam$.
We emphasize that this does not contradict the general call in \citet
{Meng2012} to follow the probabilistically-principled methods (such as
MLE and Bayes recipes) to prevent violations of self-efficiency,
precisely because the well-established principles of single-phase
inference may need to be ``re-principled'' before they can be equally
effective in the far more complicated multiphase setting.

\subsubsection*{Robustness and nuisance parameters}
In the simplest case, if a $T$ can be found such that it is a pivot
with respect to $\obsparam$ and remains dependent upon $\param$, then
sensitivity to the behavior of $\xi$ can be eliminated by preprocessing.
In such cases, an MLE or Bayes rule based on $T$ can dominate that
based on $Y$ even asymptotically.
One such example would be providing $z$-statistics from each of a set
of experiments to the downstream analyst.
This clearly limits the range of feasible downstream inferences.
With these $z$-statistics, detection of signals via multiple testing
(e.g., \citet{Benjamini1995}) would be straightforward, but efficient
combination of information across experiments could be difficult.
This is a ubiquitous trade-off of preprocessing: reductions that remove
nuisance parameters and improve robustness necessarily reduce the
amount of information available from the data.
These trade-offs must be considered carefully when designing
preprocessing techniques---universal utility is unattainable without
the original data.

A more subtle case involves the selection of $T$ as a ``partial pivot''.
In some settings, there exists a decomposition of $\obsparam$ as
$(\obsparam_1, \obsparam_2)$ such that $\dim(\obsparam_1) < D$ for some
fixed $D$ and all $N_Y$, and the distribution of $T$ is free of
$\obsparam_2$ for all values of $\obsparam_1$.
Many normalization techniques used in the microarray application of
Section~\ref{sec:examples} can be interpreted in this light.
These methods attempt to reduce the unbounded set of
experiment-specific nuisance parameters affecting $T$ to a bounded,
manageable size.

For example, suppose each processor $i$ observes $y_{ij} \sim N(\beta_0
+ \beta_{1i} x_j, \sigma^2)$, $j=1 ,\ldots,m$.
The downstream analyst wants to estimate $\beta_0$, considering $\{
\beta
_{1i} : i=1 ,\ldots,n\}$ and $\sigma^2$ as nuisance parameters.
In our previous notation, we have $\theta= \beta_0$ and $\xi=  (
\sigma^2, \beta_{11}, \ldots, \beta_{1n}  )$.
Suppose each preprocessor reduces her data to $T_i = \frac{1}{m} \sum_{j=1}^{m} (y_{ij} - \hat{\beta}_{1i} x_j)$, where $\hat{\beta}_{1i}$
is the OLS estimator of $\beta_{1i}$ based on $\{y_{ij} : j=1 ,\ldots
,m\}$.
The distribution of each $T_i$ depends on $\sigma^2$ but is free of
$\beta_{1i}$.
Hence, $T = \{T_i : i=1, \ldots,n\}$ is a partial pivot as defined
above, with $\xi_1 = \sigma^2$ and $\xi_2 = \{\beta_{1i} : i=1,
\ldots,
n\}$.

Such pivoting techniques can allow $\est(T)$ to possess favorable
properties even when $\est(Y)$ is inconsistent or grossly inefficient.
As mentioned before, this kind of careful preprocessing can dominate
Bayesian procedures in the presence of nuisance parameters when $\dim
(\obsparam)$ can grow with $N_y$.
In these regimes, informative priors on $\obsparam$ can affect
inferences even asymptotically.
However, reducing $Y$ to $T$ so only the $\obsparam_{1}$-part of
$\obsparam$ is relevant for $T$'s distribution allows information to
accumulate on $\obsparam_1$, making inferences far more robust to the
preprocessor's beliefs about $\obsparam$.

These techniques share a common conceptual framework: invariance.
Invariance has a rich history in the Bayesian literature, primarily as
a motivation for the construction of noninformative or reference priors
(e.g., \citet{Jeffreys1946}, \citet{Hartigan1964},
\citet{Geisser1979},
\citet{Berger1992},
\citet{Kass1996}).
It is fundamental to the pivotal methods discussed above and arises in
the theory of partial likelihood (\citet{Cox1975}).
We see invariance as a core principle of preprocessing, although its
application is somewhat different from most Bayesian settings.
We are interested in finding functions of the data whose distributions are invariant
to subsets of the parameter, not priors invariant to reparameterization.
For instance, the rank statistics that form the basis for Cox's
proportional hazards regression in the absence of censoring (\citeyear
{Cox1972}) can be obtained by requiring a statistic invariant to
monotone transformations of time.
Indeed, Cox's regression based on rank statistics can be viewed as an
excellent example of eliminating an infinite dimensional nuisance
parameter, i.e., the baseline hazard, via preprocessing, which retains
only the rank statistics.
The relationship between invariance in preprocessing, modeling, and
prior formulation is a rich direction for further investigation.

An interesting practical question arises from this discussion of
robustness: how realistic is it to assume efficient inference with
preprocessed data?
This may seem unrealistic as preprocessing is frequently used to
simplify problems so common methods can be applied.
However, preprocessing can make many assumptions more appropriate.
For example, aggregation can make normality assumptions more realistic,
normalization can eliminate nuisance parameters, and discretization
greatly reduces reliance on parametric distributional assumptions altogether.
It may therefore be more appropriate to assume that efficient
estimators are generally used with preprocessed data than with raw data.

The results and examples explored here show that preprocessing is a
complex topic in even large-sample settings.
It appears formally futile (but practically useful) in standard
asymptotic regimes. Under other realistic asymptotic regimes,
preprocessing emerges as a powerful tool for addressing nuisance
parameters and improving the robustness of inferences.
Having established some of the formal motivation and trade-offs for
preprocessing, we discuss further extensions of these ideas into more
difficult settings in Section~\ref{sec:future}.

\subsection{Giving all that you can}
\label{sec:completeclass}

In some cases, effective preprocessing techniques are quite apparent.
If $\pobs(Y \given X, \obsparam)$ forms an exponential family with
parameter $X$ or $(X, \obsparam)$, then we have a straightforward
procedure: retain a minimal sufficient statistic.
To be precise, we mean that one of the following factorizations holds
for a sufficient statistic $T(Y)$ of bounded dimension:
\begin{eqnarray*}
\pobs(Y \given X, \obsparam) &=& g(Y) \exp \bigl(T(Y)^\top
f(X, \obsparam ) + h(X, \obsparam) \bigr);
\\
\pobs(Y \given X, \obsparam) &=& g(Y; \obsparam) \exp
\bigl(T(Y)^\top f(X) + h(X) \bigr) .
\end{eqnarray*}
Retaining this sufficient statistic will lead to a lossless
compression, assuming that the first-phase model is correct.
Unfortunately, such nice cases are rare.
Even the Bayesian approach offers little reprieve.
Integrating $\pobs(Y \given X, \obsparam)$ with respect to a prior
$\pi
_{\obsparam}(\obsparam)$ typically removes the observation model from
the exponential family---consider, for instance, a normal model with
unknown variance becoming a $t$ distribution.

If $\log\pobs(Y \given X)$ is approximately quadratic as a function of
$X$, then retaining its mode and curvature would seem to provide much
of the information available from the data to downstream analysts.
However, such intuition can be treacherous.
If a downstream analyst is combining inferences from a set of
experiments, each of which yielded an approximately quadratic
likelihood, the individual approximations may not be enough to provide
efficient inferences.
Approximations that hold near the mode of each experiment's likelihood
need not hold away from these modes---including at the mode of the
joint likelihood from all experiments.
Thus, remainder terms can accumulate in the combination of such
approximations, degrading the final inference on $\param$.
Furthermore, the requirement that $\log\pobs(Y \given X)$ be
approximately quadratic in $X$ is quite stringent.
To justify such approximations, we must either appeal to asymptotic
results from likelihood theory or confine our attention to a narrow
class of observation models $\pobs(Y \given X)$.
Unfortunately, asymptotic theory is often an inappropriate
justification in multiphase settings, because
$X$ grows in dimension with $Y$ in many asymptotic regimes of interest,
so there is no general reason to expect information to accumulate on $X$.
These issues are of particular concern as such quadratic approximations
are a standard implicit justification for passing point estimates with
standard errors onto downstream analysts.

Moving away from these cases, solutions become less apparent.
No processing (short of passing the entire likelihood function) will
preserve all information from the sample when sufficient statistics of
bounded dimension do not exist.
However, multiphase approaches can still possess favorable properties
in such settings.

We begin by considering a stubborn downstream analyst---she has her
method and will not consider anything else.
For example, this analyst could be dead set on using linear
discriminant analysis or ANOVA.
The preprocessor has only one way to affect her results: carefully
designing a particular $T$ given to the downstream analyst.
Such a setting is extreme.
We are saying that the downstream analyst will charge ahead with a
given estimator regardless of her input with neither reflection nor judgment.
We investigate this setting because it maximizes the preprocessor's
burden in terms of her contribution to the final estimate's quality.
Formally, we consider a fixed second-stage estimator $\est(T)$; that
is, the form of its input $T$ and the function producing $\est$ are
fixed, but the mechanism actually used to generate $T$ is not.
$T$ could be, for example, a vector of fixed dimension.

As we discuss below, admissible designs for the first-phase with a
fixed second-phase method are given by a (generalized) Bayes rule.
This uses the known portion of the model $\pobs(Y \given X, \obsparam)$
to construct inputs for the second stage and assumes that any prior the
preprocessor has on $\obsparam$ is equivalent to what a downstream
analyst would have used in the preprocessor's position.
Formally, this describes all rules that are admissible among the class
of procedures using a given second-stage method, following from
previous complete class results in statistical decision theory (e.g., \citet{Berger1985},
\citet{Farrell1968}).

\subsubsection*{Admissibility}
Assume that the second-stage procedure $\est(T)$ is fixed as discussed
above and we are operating under the model (\ref{e:model}).
Further assume that the preprocessor's prior on $\obsparam$ is the only
such prior used in all Bayes rule constructions.
For $T \in\mathbb{R}^d$, consider a smooth, strictly convex loss
function~$L$.
Then, under appropriate regularity conditions (e.g., \citet{Berger1985},
\citet{Farrell1968}), if $\est(T)$ is a smooth function of
$T$, then all admissible procedures for generating $T$ are Bayes or
generalized Bayes rules with respect to the risk $R(\est(T), \param_0)$.
The same holds when $T$ is restricted to a finite set.

This guideline follows directly from conventional complete class
results in decision theory.
We omit technical details here, focusing instead on the guideline's
implications.
However, a sketch of its proof proceeds along the following lines.

There are two ways to approach this argument: intermediate loss and geometry.
The intermediate loss approach uses an intermediate loss function
$\tilde{L}(T, \param_0) = L(\est(T), \param_0)$.
This $\tilde{L}$ is the loss facing the preprocessor given a fixed
downstream procedure $\est(T)$.
If $\tilde{L}$ is well-behaved, in the sense of satisfying standard
conditions (strict convexity, or a finite parameter space, and so on),
then the proof is complete from previous results for real $T$.
Similarly, if $T$ is restricted to a finite discrete set, then we face
a classical multiple decision problem and can apply previous results to
$\tilde{L}(T, \param_0)$.
These straightforward arguments cover a wide range of realistic cases,
as \citet{Berger1985} has shown.
Otherwise, we must turn to a more intricate geometric argument.
Broadly, this construction uses a convex hull of risks generated by
attainable rules.

This guideline has direct bearing upon the development of inputs for
machine learning algorithms, typically known as \emph{feature engineering}.
Given an algorithm that uses a fixed set of inputs, it implies that
using a correctly-specified observation model to design these inputs is
necessary to obtain admissible inferences.
Thus, it is conceptually similar to ``Rao-Blackwellization'' over part
of a probability model.

However, several major caveats apply to this result.
First, on a practical level, deriving such Bayes rules is quite
difficult for most settings of interest.
Second, and more worryingly, this result's scope is actually quite limited.
As we discussed in Section~\ref{sec:missinfo}, even Bayesian estimators
can be inconsistent in realistic multiphase regimes.
However, these estimators are still admissible, as they cannot be
dominated in risk for particular values of the nuisance parameters $\xi$.
Admissibility therefore is a minimal requirement; without it, the
procedure can be improved uniformly, but with it, it can still behave
badly in many ways. Finally, there is the problem of robustness.
An optimal input for one downstream estimator $\est_1(T)$ may be a
terrible input for another estimator $\est_2(T)$, even if $\est_1$ and
$\est_2$ take the same form of inputs.
Such considerations are central to many real-world applications of
preprocessing, as researchers aim to construct databases for a broad
array of later analyses.
However, this result does show that engineering inputs for downstream
analyses using Bayesian observation models can improve overall inferences.
How to best go about this in practice is a rich area for further work.

\subsection{Counterexamples and conundrums}
\label{sec:counterexamples}

As befits first steps, we are left with a few loose ends and puzzles.
Starting with the DSC condition (\ref{eq:dsc}) of Section~\ref{sec:sufficiency}, we provide a simple counterexample to its necessity.

Suppose we have $Y_1, Y_2, X_1, X_2 \in\mathbb{R}^n$. Let $ Y_i | X_i
\sim N(X_i, I)$ independent of each other.
Now, let $X_1 = \param Z_1$, $Z_1 \sim N(0, I)$, $X_2 = \param
\mathrm{abs}(Z_2) \circ\sign(X_1)$, where $Z_2 \sim N(0, I)$, $Z_2
\indep Z_1$, $\sign(X_1)$ is a vector of signs $(-1, 0$, or $1)$ for
$X_1$, $\mathrm{abs}()$ denotes the element-wise absolute value, and
$\circ$ denotes the Hadamard product.
We fix $\param> 0$.

As our working model, we posit that $X_i |\eta\sim N(0, \workparam_i
I)$ independently.
Then, we clearly have $(Y_1^\top Y_1, Y_2^\top Y_2) = (T_1, T_2)$ as a
sufficient statistic for both $\workparam$ and $\param$.
However, the DSC does not hold for this working model.
We cannot write the actual joint distribution of $X$ as a
marginalization of $\pwork(X |\workparam)$ with respect to some
distribution over $\workparam$ in such a way that $(T_1, T_2)$ is
sufficient for $\workparam$.
To enforce $\sign(X_1) = \sign(X_2)$ under the working model, any such
model must use $\workparam$ to share this information.

For this example, we can obtain a stronger result: no factored working
model $\pwork(X |\workparam)$ exists such that (1) $Y_i^\top Y_i$ is
sufficient for $g_i(\workparam)$ under $\pywork(Y_i | g_i(\workparam))$
and (2) the DSC holds.
For contradiction, assume such a working model exists.
Under this working model, $Y_i$ is conditionally independent of
$g_i(\eta)$ given $Y_i^\top Y_i$, so we can write $\pywork(Y_i |
g_i(\eta)) = \pywork(Y_i | Y_i^\top Y_i) h_i(Y_i^\top Y_i; g_i(\eta))$.
As the DSC holds for this working model, we have
\[
\pobs(Y | \theta) = \Biggl[ \prod_{i=1}^2
\pywork\bigl(Y_i | Y_i^\top Y_i
\bigr) \Biggr] \int_{\eta} \Biggl[ \prod
_{i=1}^2 h_i\bigl(Y_i^\top
Y_i; g_i(\eta)\bigr) \Biggr] p_\workparam(\dd
\workparam| \param) .
\]
Hence, we must have $Y_1$ conditionally independent of $Y_2$ given
$(Y_1^\top Y_1, Y_2^\top Y_2)$.
However, this conditional independence does not hold under the true model.
Hence, the given working model cannot both satisfy the DSC and have
$Y_i^\top Y_i$ sufficient for each $g_i(\eta)$.

The issue here is unparameterized dependence, as mentioned in Section~\ref{sec:sufficiency}.
The $X$'s have a dependence structure that is not captured by $\param$.
Thus, requiring that a working model preserves sufficiency for $\param$
does not ensure that it has enough flexibility to capture the true
distribution of $Y$.
A weaker condition than the DSC (\ref{eq:dsc}) that is necessary and
sufficient to ensure that all sufficient statistics for $\workparam$
are sufficient for $\param$ may be possible.

From Sections~\ref{sec:missinfo} and~\ref{sec:completeclass}, we are
left with puzzles rather than counterexamples.
As mentioned previously, many optimality results are trivial without
sufficient constraints.
For instance, minimizing risk or maximizing Fisher information naively
yield uninteresting (and impractical) multiphase strategies: have the
preprocessor compute optimal estimators, then pass them downstream.
Overly tight constraints bring their own issues.
Restricting downstream procedures to excessively narrow classes (e.g.,
point estimates with standard errors) limits the applied utility of
resulting theory and yields little insight on the overall landscape of
multiphase inference.
Striking the correct balance with these constraints is a core challenge
for the theory of multiphase inference and will require a combination
of computational, engineering, and statistical insights.

\section{From the past to the future}
\label{sec:remarks}

As we discussed in Sections~\ref{sec:concepts} and~\ref{sec:theory}, we
have a deep well of questions that motivate further research on
multiphase inference.
These range from the extremely applied (e.g., enhancing preprocessing
in astrophysical systems) to the deeply theoretical (e.g., bounding
the performance of multiphase procedures in the presence of nuisance parameters and
computational constraints).
We outline a few directions for this research below.

But, before we look forward, we take a moment to look back and place
multiphase inference within the context of broader historical debates.
Such ``navel gazing'' helps us to understand the connections and
implications of the theory of multiphase inference.

\subsection{Historical context}

On a historical note, the study of multiphase inference touches the
long-running debate over the role of decision theory in statistics.
One side of this debate, championed by Wald and Lehmann (among others),
has argued that decision theory lies at the core of statistical inference.
Risk-minimizing estimators and, more generally, optimal decision rules
play a central role in their narrative.
Even subjectivists such as Savage and de Finetti have embraced the
decision theoretic formulation to a large extent.
Other eminent statisticians have objected to such a focus on decisions.
As noted by \citet{Savage1976}, Fisher in particular vehemently
rejected the decision theoretic formulation of statistical inference.
One interpretation of Fisher's objections is that he considered
decision theory useful for eventual economic decision-making, but not
for the growth of scientific knowledge.

We believe that the study of multiphase inference brings a unifying
perspective to this debate.
Fisher's distinction between intermediate processing and final
decisions is fundamental to the problem of multiphase inference.
However, we also view decision theory as a vital theoretical tool for
the study of multiphase inference.
Passing only risk-minimizing point estimators to later analysts is
clearly not a recipe for valid inference.
The key is to consider the use of previously generated results
explicitly in the final decision problem.
In the study of multiphase inference, we do so by focusing on the
separation of knowledge and objectives between agents.
Such separation between preprocessing and downstream inference maps
nicely to Fisher's distinction between building scientific knowledge
and reaching actionable decisions.

Thus, we interpret Fisher's line of objections to decision-theoretic
statistics as, in part, a rejection of adopting a myopic single-phase
perspective in multiphase settings.
We certainly do not believe that our work will bring closure to such an
intense historical debate.
However, we do see multiphase inference as an important bridge between
these competing schools of thought.

\subsection{Where can multiphase inference go from here?}
\label{sec:future}

We see a wide range of open questions in multiphase inference.
Can more systematic ways to leverage the potential of preprocessing be
developed?
Is it possible to create a mathematical ``warning system,'' alerting
practitioners when their inferences from preprocessed data are subject
to severe degradation and showing where additional forms of
preprocessing are required?
And, can multiphase inference inform developments in distributed
statistical computation and massive-data inference (as outlined below
in Section~\ref{sec:computation})?
All of these problems call for a shared collection of statistical
principles, theory, and methods.
Below, we outline a few directions for the development of these tools
for multiphase inference.


\subsubsection*{Passing information}
The mechanics of passing information between phases constitute a major
direction for further research.
One approach leverages the fact that the likelihood function itself is
always a minimal sufficient statistic.
Thus, a set of (computationally) efficient approximations to the
likelihood function $L(X, \obsparam; Y)$ for $(X, \xi)$ could provide
the foundation for a wide range of multiphase methods.
Many probabilistic inference techniques for the downstream model (e.g.,
MCMC samplers) would be quite straightforward to use given such an
approximation.
The study of such multiphase approximations also offers great dividends
for distributed statistical computation, as discussed below.
We believe these approximations are promising direction for
general-purpose preprocessing.
However, there are stumbling blocks.

First, nuisance parameters remain an issue.
We want to harness and understand the robustness benefits offered by
preprocessing, but likelihood techniques themselves offer little
guidance in this direction.
Even the work of \citet{Cox1975} on partial likelihood focuses on the
details of estimation once the likelihood has been partitioned.
We would like to identify the set of formal principles underlying
techniques such as partial pivoting (to mute the effect of
infinite-dimensional nuisance parameters), building a more rigorous
understanding of the role of preprocessing in providing robust inferences.
As discussed in Section~\ref{sec:missinfo}, invariance relationships
may be a useful focus for such investigations, guiding both Bayesian
and algorithmic developments.

Second, we must consider the burden placed on downstream analysts by
our choice of approximation.
Probabilistic, model-based techniques can integrate such information
with little additional development.
However, it would be difficult for a downstream analyst accustomed to,
say, standard regression methods to make use of a complex emulator for
the likelihood function.
The burden may be substantial for even sophisticated analysts.
For instance, it could require a significant amount of effort and
computational sophistication to obtain estimates of $X$ from such an
approximation, and estimates of $X$ are often of interest to downstream
analysts in addition to estimates of $\param$.

\subsubsection*{Bounding errors and trade-offs}
With these trade-offs in mind and through the formal analysis of
widely-applicable multiphase techniques, we can begin to establish
bounds on the error properties of such techniques in a broad range of
problems under realistic constraints (in both technical and human terms).
More general constraints, for instance, can take the form of upper
bounds on the regret attainable with a fixed amount of information
passed from preprocessor to downstream analyst for fixed classes of
scientific models.
Extensions to nonparametric downstream methods would have both
practical and theoretical implications.
In cases where the observation model is well-specified but the
scientific model is less clearly defined, multiphase techniques can
provide a useful alternative to computationally-expensive
semi-parametric techniques.
Fusing principled preprocessing with flexible downstream inference may
provide an interesting way to incorporate model-based subject-matter
knowledge while effectively managing the bias-variance trade-off.

\subsubsection*{Links to multiple imputation}
The directions discussed above share a conceptual, if not technical,
history with the development of congeniality (\citet{Meng1994}).
Both the study of congeniality in MI and our study of multiphase
inference seek to bound and measure the amount of degradation in
inferences that can occur when agents attempt (imperfectly) to combine
information.
Despite these similarities, the treatment of nuisance parameters are
rather different.
Nuisance parameters lie at the very heart of multiphase inference,
defining many of its core issues and techniques.
For MI, the typical approaches have been to integrate them out in a
Bayesian analysis (e.g., \citet{Rubin1996}) or assume that the final
analyst will handle them (e.g., \citet{Nielsen2003}).
Recent work by \citet{Xie2012} has shed new light on the role of
nuisance parameters in MI, but the results are largely negative,
demonstrating that nuisance parameters are often a stumbling block for
practical MI inference.
Understanding the role of preprocessing in addressing nuisance
parameters, providing robust analyses, and effectively distributing
statistical inference represent further challenges beyond those pursued
with MI.
Therefore, much remains to be done in the study of multiphase
inference, both theoretical and methodological.

\subsection{How does multiphase inference inform computation?}
\label{sec:computation}

We also see multiphase inference as a source for computational
techniques, drawing inspiration from the history of MI.
MI was initially developed as a strategy for handling missing data in
public data releases.
However, because MI separates the task of dealing with incomplete data
from the task of making inferences, its use spread.
It has frequently been used as a practical tool for dealing with
missing-data problems where the joint inference of missing data and
model parameters would impose excessive modeling or computational burdens.
That is, increasingly the MI inference is carried out from imputation
through analysis by a single analyst or research group.
This is feasible as a computational strategy only because the error
properties and conditions necessary for the validity of MI are
relatively well-understood (e.g., \citet{Meng1994},
\citet{Xie2012}).\looseness=1

Multiphase methods can similarly guide the development of efficient,
statistically-valid computational strategies.
Once we have a theory showing the trade-offs and pitfalls of multiphase
methods, we will be equipped to develop them into general computational
techniques.
In particular, our experience suggests that models with a high degree
of conditional independence (e.g., exchangeable distributions for $X$)
can often provide useful inputs for multiphase inferences, even when
the true overall model has a greater degree of stochastic structure.
The conditional independence structure of such models allows for highly
parallel computation with first-phase procedures, providing huge
computational gains on modern distributed systems compared to methods
based on the joint model.\looseness=1

For example, in \citet{Blocker2012}, a factored model was used to
preprocess a massive collection of irregularly-sampled astronomical
time series.
The model was sophisticated enough to account for complex observation
noise, yet its independence structure allowed for efficient
parallelization of the necessary computation.
Its output was then combined and used for population-level analyses.
Just as Markov chain Monte-Carlo (MCMC) has produced a windfall of
tools for approximate high-dimensional integration (see \citet{Brooks2010} for many
examples), we believe that this type of principled
preprocessing, with further theoretical underpinnings, has the
potential to become a core tool for the statistical analysis of massive
datasets.\looseness=1

\section*{Acknowledgments}
We would like to acknowledge support from the Arthur P. Dempster Award
and partial financial support from the NSF.
We would also like to thank Arthur P. Dempster and Stephen Blyth for
their generous feedback.
This work developed from the inaugural winning submission for said award.
We also thank David van Dyk, Brandon Kelly, Nathan Stein, Alex D'Amour,
and Edo Airoldi for valuable discussions and feedback, and Steven Finch
for proofreading.
Finally, we would like to thank our reviewers for their thorough and
thoughtful comments, which have significantly enhanced this
paper.\looseness=1

%



\printhistory

\end{document}